\documentclass[9pt]{article}
\usepackage[utf8]{inputenc}
\usepackage[english]{babel}

\usepackage{eso-pic}
\usepackage{amsmath} 
\usepackage{amssymb}
\usepackage{mathrsfs}
\usepackage{amsthm}
\usepackage{remreset}
\usepackage{tikz}
\usetikzlibrary{calc} 
\usetikzlibrary{math}
\usepackage{tikz-cd} 
\usepackage{graphicx}
\usetikzlibrary{shapes.misc}
\usepackage{appendix}
\usepackage{xcolor}
\usepackage{graphicx}
\usepackage{lineno}
\usepackage{titlesec}
\usepackage{scalerel}[2016/12/29]
\usepackage[numbers]{natbib}
\usepackage{hyperref}
\usepackage{subfiles}
\usepackage[all,cmtip]{xy}

\usepackage{setspace}

\newcommand{\poubelle}[1]{}

\newcommand{\N}{\mathbb{N}}

\newcommand{\U}{\mathcal{U}}

\newcommand{\A}{\mathcal{A}}

\newcommand{\F}{\mathbb{F}_2}

\newcommand{\Fc}{\mathcal{F}}
\newcommand{\Fo}{\mathcal{F}_\omega}
\newcommand{\K}{\mathcal{K}}
\newcommand{\Fk}[1]{\mathcal{F}^{<#1}}
\newcommand{\Fok}[1]{\mathcal{F}_\omega^{<#1}}

\newcommand{\E}{\mathcal{V}^f}
\newcommand{\Ef}{\mathcal{V}}

\newcommand{\Hom}{\text{Hom}}

\newcommand{\HomK}{\text{Hom}_\K}
\newcommand{\HomF}{\text{Hom}_{\F}}
\newcommand{\HomU}{\text{Hom}_\U}

\newcommand{\Sq}{\text{Sq}}
\newcommand{\Nil}{\mathcal{N}il}

\newcommand{\id}{\text{id}}

\newcommand{\Ima}{\text{Im}}

\newcommand{\coker}{\text{coker}}

\newcommand{\PFVF}{\mathcal{P}\text{fin}^{(\mathcal{V}^f)^\text{op}}}
\newcommand{\Skop}[1]{(\E)^{\mathcal{S}(#1)^{op}}}

\def\commutatif{\ar@{}[rd]|{\circlearrowleft}}

\newtheorem{theorem}{Theorem}[section]

\newtheorem{proposition}[theorem]{Proposition}
\newtheorem{lemma}[theorem]{Lemma}
\newtheorem{corollary}[theorem]{Corollary}
\newtheorem{innercustomthm}{Theorem}
\newenvironment{customthm}[1]
{\renewcommand\theinnercustomthm{#1}\innercustomthm}
{\endinnercustomthm}
\newtheorem{ppp}{Proposition}
\newenvironment{pp}[1]
{\renewcommand\theppp{#1}\ppp}
{\endppp}

\theoremstyle{definition}

\newtheorem{definition}[theorem]{Definition}
\newtheorem{remark}[theorem]{Remark}

\theoremstyle{definition}
\newtheorem{example}[theorem]{Example}

\theoremstyle{definition}
\newtheorem{notation}[theorem]{Notation}

\usepackage{natbib}
\usepackage{graphicx}
\usepackage[all]{xy}
\usepackage{xypic}
\usepackage{amssymb}
\usepackage{amsmath}
\usepackage{mathrsfs}
\usepackage{mathtools}
\usepackage{amsthm}
\usepackage{moreverb}
\usepackage{geometry}

\title{APPROXIMATION OF THE CENTRE OF UNSTABLE ALGEBRAS USING THE NILPOTENT FILTRATION}
\author{OURIEL BLOEDE}

\date{}
\begin{document}

\maketitle

\begin{abstract}
In \cite{B1}, we computed the set $\textbf{C}(K)$ of central elements of an unstable algebra $K$ over the Steenrod algebra, in the sense of Dwyer and Wilkerson, when $K$ is noetherian and $nil_1$-closed.

For $K$ noetherian and $k$ a positive integer, we define $\textbf{C}_k(K)$, the set of so-called central elements of $K$ away from $\Nil_k$ in such a way that, for $K$ $nil_k$-closed, $\textbf{C}(K)=\textbf{C}_k(K)$.

The sets $\textbf{C}_k(K)$ are a decreasing filtration, and we describe the obstruction for an element in $\textbf{C}_k(K)$ to be in $\textbf{C}_{k+1}(K)$. Since, for $K$ noetherian, $K$ is always $nil_k$-closed for $k$ big enough, this gives us a way to compute the set of central elements of $K$.
\end{abstract}

\section{Introduction}

For $\A$ the Steenrod algebra modulo $2$, $\U$ and $\K$ denote the categories of unstable modules and unstable algebras over $\A$.\\

In \cite{DW1}, Dwyer and Wilkerson defined the centre of an unstable algebra over the Steenrod algebra. In the case where $K$ is noetherian and connected, the set of central elements of $K$ coincides with the set of pairs $(V,\phi)$ such that
\begin{enumerate}
    \item $\phi\in\HomK(K,H^*(V))$,
    \item $K$ admits a structure $\kappa$ of $H^*(V)$-comodule in $\K$, such that the following diagram commutes:
    $$\xymatrix{ 
K\ar[rr]^\kappa\ar[rrd]_\phi & & K\otimes H^*(V)\ar[d]^-{\epsilon_K\otimes\id} & \\
& & H^*(V),}$$ where $\epsilon_K$ denotes the augmentation of $K$ (which is uniquely defined because of the connectedness of $K).$
\end{enumerate}
For $X$ a topological space they used the centre of $H^*(X)$ (the mod $2$ singular cohomology of $X$) to prove an analogue of Sullivan's conjecture, proved by Miller in \cite{M1}. Namely, they proved that, for $V$ an elementary abelian $2$-group, $BV$ its classifying space, $X$ $1$-connected and $2$-complete and $\gamma$ a continuous map from $BV$ to $X$, the connected component of $\gamma$ in the mapping space $\text{Map}_*(BV,X)$ is weakly contractible if and only if $(V,\gamma^*)$ is a central element of $H^*(X)$, for $\gamma^*$ the morphism induced by $\gamma$ from $H^*(X)$ to $H^*(V):=H^*(BV)$.\\

In \cite{K11} (for the cohomology of groups) and in \cite{heard2019depth} (for any unstable algebra), Kuhn and Heard used central elements of an unstable algebra to find bounds to their depth.\\

In a similar way, for $K$ an unstable algebra, one can define $M$-central elements for $M$ in $K-\U$ the category of modules over $K$ in $\U$ (see \cite{heard2019depth}). We will denote by $\textbf{C}(K)$, respectively $\textbf{C}(M;K)$, the set of pairs $(V,\phi)$ with $V$ a finite dimensional $\F$-vector space and $\phi\ :\ K\rightarrow H^*(V)$ a morphism of unstable algebras, such that $(V,\phi)$ is a central element, respectively a $M$-central element.\\

The category $\U$ of unstable algebras over the Steenrod algebra admits a filtration by localizing subcategories $\U=\Nil_0\supset\Nil_1\supset ...\supset\Nil_k\supset ...$ (see \cite{S1}). In \cite{HLS2} and \cite{HLS1}, Henn, Lannes and Schwartz constructed category equivalences between the categories $\U/\Nil_k$ and categories of functors denoted by $\Fok{k}$. In \cite{B1}, we used the equivalence between the category $\U/\Nil_1$ and the category $\Fok{1}$ to compute the set of central elements of $nil_1$-closed, noetherian, integral, unstable algebras. In this article we want to show how to approximate the sets of central elements $\textbf{C}(K)$ and $\textbf{C}(M;K)$, for $K$ non $nil_1$-closed, using the nilpotent filtration.\\

We define a notion of central elements away from $\Nil_k$ and denote by $\textbf{C}_k(K)$ and $\textbf{C}_k(M;K)$ the sets of central elements and $M$-central elements away from $\Nil_k$. For $K$ an unstable algebra and $M\in K-\U$, we prove that $l_k(K)$ (the $nil_k$-closure of $K$) and $l_k(M)$ (the $nil_k$-closure of $M$) are respectively an unstable algebra and an object of $l_k(K)-\U$ (see Theorem \ref{princ1} and Proposition \ref{princ2}). Since for all finite dimensional vector space $V$, $H^*(V)$ is both injective and $nil_1$-closed, every morphism of unstable algebras $K\rightarrow H^*(V)$ factorises uniquely through $l_k(K)$. Under this identification we find that the sets of central and $M$-central elements away from $\Nil_k$ are $\textbf{C}(l_k(K))$ and $\textbf{C}(l_k(M);l_k(K))$ (Propositions \ref{nilfer} and \ref{nilfer2}).\\

By construction, we have, for all $k\in\N$, the inclusions $\textbf{C}_{k+1}(K)\subset\textbf{C}_k(K)$ and $\textbf{C}_{k+1}(M;K)\subset\textbf{C}_k(M;K)$. This filtration gives us a good way to approximate $\textbf{C}(K)$ and $\textbf{C}(M;K)$. 

\begin{pp}{\ref{princ3}}
Let $K\in\K$ and $M\in K-\U$, then $$\textbf{C}(K)=\bigcap\limits_{k\in\N^*}\textbf{C}_k(K)$$
and $$\textbf{C}(M;K)=\bigcap\limits_{k\in\N^*}\textbf{C}_k(M;K).$$
\end{pp}

The main theorems of this article describe the obstruction for an element in $\textbf{C}_k(K)$ (respectively in $\textbf{C}_k(M;K)$) to be in $\textbf{C}_{k+1}(K)$ (respectively in $\textbf{C}_{k+1}(M;K)$). We consider the following exact sequence, which is an exact sequence in $K-\U$.
$$
\xymatrix{0\ar[r]  &\ker(\lambda_k)\ar[r]^-{\lambda_k} & l_{k+1}(K)\ar[r] & l_k(K)\ar[r]  & \coker(\lambda_k)\ar[r] & 0.}$$

\begin{customthm}{\ref{gaga}}
Let $K$ be an unstable algebra and let $\phi\in\HomK(K,H^*(W))$. Then, the two following conditions are equivalent:
\begin{enumerate}
    \item $(W,\phi)\in\textbf{C}_{k+1}(K)$,
    \item $(W,\phi)\in\textbf{C}_k(K)$, $(W,\phi)\in\textbf{C}_{k+1}(\ker(\lambda_k);K)$ and  $(W,\phi)\in\textbf{C}_{k+1}(\coker(\lambda_k);K)$.
\end{enumerate}
\end{customthm}

We have a similar statement for $M\in K-\U$.\\

The interest of this characterisation depends on the fact that $\ker(\lambda_k)$ and $\coker(\lambda_k)$ are both objects of $\Nil_k$. In the fifth section, we explain using the category equivalences of Henn, Lannes and Schwartz, how to compute $\textbf{C}_1(K)$ (this was already studied extensively in \cite{B1}) as well as $\textbf{C}_{k+1}(M;K)$ for $M$ a $k$-nilpotent object in $K-\U$. For $K$ noetherian, using that there is some integer $t$ such that $K$ is $nil_t$-closed, Theorem \ref{gaga} gives us an algorithm that computes $\textbf{C}(K)$.\\

In the last section, we give various example of computations of $\textbf{C}(K)$ using Theorem \ref{gaga}.\\

\textbf{Acknowledgements: }I want to thank Geoffrey Powell for his careful proofreading. This work was partially supported by the ANR Project {\em ChroK}, {\tt ANR-16-CE40-0003}.

\section{Unstable modules over the Steenrod algebra and the Nilpotent filtration}

In this section, we recall some known facts about Lannes' $T$ functor as well as results from \cite{HLS1} about the localization of the category of unstable modules away from $k$-nilpotent objects. Recollections about unstable algebras, unstable modules and $k$-nilpotent objects can be found in \cite{S1}. In the following, $\A$ denotes the Steenrod algebra over $\F$, $\U$ and $\K$ denote the category of unstable modules and unstable algebras over $\A$ and $\Nil_k$ denotes the class of $k$-nilpotent objects in $\U$.

\subsection{Standard Projective objects in $\U$}

We recall the definition of $F(n)$ from \cite{S1}. 

\begin{proposition}\cite[Proposition 1.6.1]{S1}\label{projst}
There is, up to isomorphism, a unique unstable module $F(n)$ with a class $x(n)$ of degree $n$ such that the natural transformation $f\mapsto f(x(n))$ from $\HomU((F(n),M)$ to $M^n$ is a natural isomorphism.
\end{proposition}

$F(n)$ is projective, since the functor which maps $M$ to $M^n$ is exact. Furthermore, $F(n)$ is generated by $x(n)$ as a $\A$-module.

\begin{definition}\label{projdef}
For $n$ and $i\leq n$ in $\N$, and for $\theta\in\A$, let $u_\theta\ :\ F(n+|\theta|)\rightarrow F(n)$ be the morphism which maps $x(n+|\theta|)$ to $\theta x(n)$.
\end{definition}

We introduce a graded unstable module $F(*)$, which will be of great use later.

\begin{definition}
Let $F(*)$ be the graded unstable module whose component of degree $n$ is $F(n)$.
\end{definition}

From Proposition \ref{projst}, we get the following isomorphism of graded vector spaces $\HomU(F(*),M)\cong M$, where the graduation of $\HomU(F(*),M)$ comes from the one of $F(*)$. The maps $u_\theta$ defined in Definition \ref{projdef} induce a structure of left $\A$-module on $F(*)$. Therefore, $\HomU(F(*),M)$ gets a structure of right $\A$-module. The following Proposition follows directly from Definition \ref{projdef}.

\begin{proposition}\label{fufu}
The isomorphism $\HomU(F(*),M)\cong M$ is an isomorphism of right $\A$-modules.
\end{proposition}

\subsection{Injective objects in $\U$}

We recall the description of injective objects in $\U$ from \cite{ASENS_1986_4_19_2_303_0}. The proofs can be found in \cite{S1}.

\begin{proposition}\cite[Part 2.2]{S1}
The functor from $\U$ to the category of vector spaces defined by $M\mapsto (M^n)^\sharp:=\Hom_{\E}(M^n,\F)$ is representable.
\end{proposition}

\begin{definition}
For $n\in\N$ let $J(n)\in\U$ be the unstable module satisfying that $\Hom_\U(\_,J(n))\cong ((\_)^n)^\sharp$. 
\end{definition}

\begin{remark}
By definition, $\Hom_\U(\_,J(n))$ is right-exact, hence $J(n)$ is injective.
\end{remark}

\begin{proposition}\cite[Theorem 3.1.1]{S1}
For every $V\in\E$ and every $n\in\N$, $H^*(V)\otimes J(n)$ is injective in $\U$.
\end{proposition}

\begin{remark}
In particular, since $J(0)$ is isomorphic as an unstable module to $\F$ concentrated in degree $0$, the latter proposition implies that $H^*(V)$ is injective.
\end{remark}

\begin{theorem}\cite[Théorème 3.11.1]{S1}\label{intro1}
Let $\mathcal{L}$ be a set containing exactly one element in each isomorphism class of indecomposable factor of $H^*(\F^n)$ as an object of $\U$, with $n$ running through $\N$. Then, for all injective object $I$ in $\U$, there is a unique family of cardinals $(a_{L,i})_{(L,i)\in\mathcal{L}\times \N}$ such that $I\cong\bigoplus\limits_{(L,i)}(L\otimes J(i))^{\oplus a_{L,i}}$.

\end{theorem}

Since the category $\U$ has enough injective objects, we have the following corollary.

\begin{corollary}\label{injreso}
Every object $M$ in $\U$ admits an injective resolution $$0\rightarrow M\rightarrow I_0\rightarrow I_1\rightarrow ...\ ,$$ such that, for each $k$, $I_k$ is a direct sum of modules of the form $H^*(V)\otimes J(i)$ with $V\in\E$ and $i\in\N$.
\end{corollary}

\subsection{The $T$ functor}

Let us recall the definition of Lannes' $T$ functor.

\begin{theorem}\cite[Proposition 2.1]{L1}
For $V$ a finite dimensional vector space, the functor $- \otimes H^*(V)$ has a left adjoint $T_V$.
\end{theorem}

\begin{proposition}\cite[Theorems 3.2.2 and 3.5.1]{S1}\label{muche}
\begin{enumerate}
    \item The functor $T_V$ is exact.
    \item For $M$ and $N$ in $\U$, $T_V(M\otimes N)$ is naturally isomorphic to $T_V(M)\otimes T_V(N)$.
\end{enumerate}
\end{proposition}

\begin{proposition}\cite[Proposition 3.8.4]{S1}\label{Tvgebre}
\begin{enumerate}
    \item Let $K$ be an unstable algebra, then $T_V(K)$ has a natural unstable algebra structure.
    \item Then, $T_V$ defines a functor from $\K$ to $\K$ which is left adjoint to the tensor product with $H^*(V)$ in $\K$.
\end{enumerate}

\end{proposition}

\begin{example}\cite[Proposition 3.3.6]{S1}
For $M\in \U$ bounded in degree, $T_V(M)\cong M$, in particular, for all $n\in\N$, $T_V(J(n))\cong J(n)$.
\end{example}

\begin{example}\cite[3.9.1]{S1}\label{hetv}
For $V$ and $W$ two finite-dimensional $\F$-vector spaces, there is an isomorphism of unstable algebras, natural in both $V$ and $W$, $T_V(H^*(W))\cong H^*(W)\otimes\F^{\Hom(V,W)}$. By the adjunction property, we get that $(\F^{\Hom(V,W)})^\sharp\cong\HomU(H^*(W),H^*(V)).$ In other words, $\F\left[\Hom(V,W)\right]\cong\HomU(H^*(W),H^*(V))$, which is a theorem first proved by Adams, Gunawardena and Miller.

\end{example}

Since $H^*(V)$ is functorial in $V$ and $H^*(V\oplus V)\cong H^*(V)\otimes H^*(V)$, the codiagonal $\nabla\ :\ V\oplus V\rightarrow V$ induces a natural coalgebra structure $\nabla^*\ :\ H^*(V)\rightarrow H^*(V)\otimes H^*(V)$ on $H^*(V)$.

\begin{definition}
 For $M\in\U$ and $V\in\E$ (the category of finite-dimensional $\F$-vector spaces), let $\kappa_{M,V}\ :\ T_V(M)\rightarrow T_V(M)\otimes H^*(V)$ be the adjoint of the composition: $$\xymatrix{M\ar[r]& T_V(M)\otimes H^*(V)\ar[rr]^-{id_{T_V(M)}\otimes \nabla^*}& & T_V(M)\otimes H^*(V)\otimes H^*(V),}$$ where the first map is the adjoint of the identity of $T_V(M)$.
\end{definition}

\begin{proposition}\cite[1.13]{HLS1}\label{conatur}
\begin{enumerate}
    \item $\kappa_{M,V}$ endows $T_V(M)$ with a $H^*(V)$-comodule structure which is natural with respect to $M$.
    \item For every $\alpha : V\rightarrow W$ in $\E$, the following diagram commutes:  
$$\xymatrix{
T_V(M) \ar[rr]^{\kappa_{M,V}}\ar[dd]_{\alpha_*}& & H^{*}(V)\otimes T_V(M) \ar[rd]^{id\otimes \alpha_*} & \\
 & & & H^{*}(V)\otimes T_W(M)\\
T_W(M) \ar[rr]^{\kappa_{M,W}}& & H^{*}(W)\otimes T_W(M) \ar[ru]_{\alpha^*\otimes id} & 
.}$$ 
\end{enumerate}

\end{proposition}

\subsection{$\Nil_k$-localisation of unstable modules}

 The class $\Nil_k$ of $k$-nilpotent modules is a Serre class in $\U$, we recall from \cite{HLS1} the existence of an equivalence of categories between $\U/\Nil_k$ (defined as in \cite{G1}) and a category of functors.\\
 
 The following is essentially due to Gabriel in \cite{G1}.

\begin{theorem}\cite[Part I.1]{HLS1}
There is an adjunction:
$$\xymatrix{r_k : \U\ar@<2pt>[r] & \ar@<2pt>[l] \U/\Nil_k: s_k,}$$
where $r_k$ is exact and such that, for $\phi$ a morphism of unstable modules, $r_k(\phi)$ is an isomorphism if and only if $\ker(\phi)$ and $\coker(\phi)$ are objects in $\Nil_k$. 

Then, $\U/\Nil_k$ satisfies the following universal property: for $A$ an abelian category and $F\ :\ \U\rightarrow A$ an exact functor such that for all $M\in\Nil_k$, $F(M)=0$, there exists a unique exact functor $G\ :\ \U/\Nil_k\rightarrow A$ such that $F=G\circ r_k$.
\end{theorem}

\begin{definition}
For $M$ an unstable module, 
\begin{enumerate}
    \item $l_k(M):=s_k\circ r_k(M)$ is the $nil_k$-localisation of $M$,
    \item $M$ is $nil_k$-closed if the unit of the adjunction $\lambda_k\ :\ M\rightarrow l_k(M)$ is an isomorphism,
    \item $M$ is $nil_k$-reduced if the unit of the adjunction $\lambda_k\ :\ M\rightarrow l_k(M)$ is injective,
    \item $nil_k(M)$ is the kernel of $\lambda_k$.
\end{enumerate}

\end{definition}

 \begin{remark}
 $nil_k(M)$ is the largest $k$-nilpotent sub-module of $M$.
 \end{remark}

\begin{definition}
\begin{enumerate}
    \item Let $\U^{<k}$ be the full subcategory of $\U$ of unstable modules concentrated in degrees less than $k$.
    \item Let $(\_)^{<k}\ :\ \U\rightarrow \U^{<k}$ be the functor which maps an unstable module $M$, to $M^{<k}$ which is the quotient of $M$ by its elements of degree greater than or equal to $k$.
\end{enumerate}
\end{definition}

\begin{remark}
 $(\_)^{<k}$ is left adjoint to the forgetful functor from $\U^{<k}$ to $\U$.
\end{remark}

\begin{proposition}
For $M$ and $N$ in $\U^{<k}$, $M\otimes_{\U^{<k}}N:=(M\otimes N)^{<k}$ endows $\U^{<k}$ with a symmetric monoidal structure.
\end{proposition}

\begin{corollary}\label{comod}
$(\kappa_{M,V})^{<k}$ is a $H^*(V)^{<k}$-comodule structure in $\U^{<k}$ on $T_V(M)^{<k}$. 
\end{corollary}

Let us recall the definition of the categories $\Fk{k}$ from \cite{HLS1}. In the following we will use the notation $H^{*<k}(V):=(H^*(V))^{<k}$.

\begin{definition}\label{fk}

Let $\mathcal{F}^{<k}$ be the category whose objects are functors $F$ from $\E$ to $\mathcal{U}^{<k}$, such that for all $V\in\E$, $F(V)$ is provided with a $H^{*<k}(V)$-comodule structure satisfying that for all $\alpha : V\rightarrow W$ the following diagram commutes: 
$$\xymatrix{
F(V) \ar[rr]^-{\kappa_{F,V}}\ar[dd]_{F(\alpha)}& & H^{*<k}(V)\otimes_{\U^{<k}} F(V) \ar[rd]^{id\otimes F(\alpha)} & \\
 & & & H^{*<k}(V)\otimes_{\U^{<k}} F(W)\\
 F(W) \ar[rr]^-{\kappa_{F,W}}& & H^{*<k}(W)\otimes_{\U^{<k}} F(W) \ar[ru]_{\alpha^*\otimes id}, & 
}$$ and whose morphisms are natural transformations $\phi : F\rightarrow G$, such that for all $V\in\E$, $\phi_V$ is a morphism of $H^{*<k}(V)$-comodules.
\end{definition}

\begin{example}
By Proposition \ref{conatur} and Corollary \ref{comod}, for $M\in\U$, $V\mapsto T_V(M)^{<k}$ defines an object in $\Fk{k}$. 
\end{example}

\begin{definition}
Let $f^{<k}\ :\ \U\rightarrow \Fk{k}$ be the functor which maps the unstable module $M$ to the functor $V\mapsto T_V(M)^{<k}$ in $\Fk{k}$.

\end{definition}

\begin{lemma}
\begin{enumerate}
    \item The category $\Fk{k}$ is abelian.
    \item Furthermore, since the functor $T_V$ is exact for every $V\in\E$, the functor $f^{<k}$ is exact.

\end{enumerate}
 
\end{lemma}

\begin{lemma}
\begin{enumerate}
    \item The structure of symmetric monoidal category on $\U^{<k}$ induces one on $\Fk{k}$.
    \item Furthermore, since the functor $T_V$ commutes with tensor products, the functor $f^{<k}$ is a strong monoidal functor.

\end{enumerate}
\end{lemma}

The class $\Nil_k$ satisfies that $f^{<k}(M)=0$ if and only if $M\in\Nil_k$. Thus, $f^{<k}$ induces a functor $\Tilde{f}^{<k}$ from $\U/\Nil_k$ to $\Fk{k}$ such that $f^{<k}=\Tilde{f}^{<k}\circ r_k$.\\

We will denote by $\Fok{k}$ the essential image of $f^{<k}$ in $\Fk{k}$.

\begin{theorem}\cite[Theorem 2.1]{HLS1}\label{equivk}
The functor $f^{<k}$ induces an equivalence of categories:
$$\xymatrix{
\Tilde{f}^{<k} :  \mathcal{U}/Nil_k\ar@<2pt>[r] & \mathcal{F}_\omega^{<k}\ar@<2pt>[l] \ : \Tilde{m}^{<k}}.$$
\end{theorem}

\begin{definition}\label{mmm}

Let $m^{<k}\ :\ \Fok{k}\rightarrow \U$ be the composition of $\Tilde{m}^{<k}$ with $s_k$.
\end{definition}

\begin{remark}
By construction, $m^{<k}$ is right adjoint to $f^{<k}$.
\end{remark}

\begin{remark}
For $k=1$, $H^{*<1}(V)\cong\F$ for all $V\in\E$. Therefore, the category $\Fc:=\Fk{1}$ is easier to describe, it is the category of functors between the category $\E$ and the category $\Ef$ of $\F$-vector spaces (not necessarily of finite dimension). In this case, $\Fo:=\Fok{1}$ is the category of analytic functors from $\E$ to $\Ef$ and $f^{<1}$ is the functor which maps $M\in\U$ to the functor which sends $V\in\E$ to $T_V(M)^0\in\Ef$.
\end{remark}

For $k\in\N^*$ and $M\in\Nil_k$, $f^{<k+1}(M)$ is concentrated in degree $k$ (since $T_V(M)^{<k}$ is trivial for any $V\in\E$). Therefore, for $V\in\E$, the $H^{*<k+1}(V)$-comodule structure on $f^{<k+1}(M)(V)$ needs to be trivial, hence $f^{<k+1}(M)$ can be reduced to the functor $f^k(M)$ which maps $V$ to $T_V(M)^k$ which is an object of $\Fo$. This gives rise to the following. 

\begin{theorem}\cite[]{S1}\label{equivknil}
The functor $f^{k}$ induces an equivalence of categories:
$$\xymatrix{
\Tilde{f}^{k} :  \Nil_{k}/\Nil_{k+1}\ar@<2pt>[r] & \mathcal{F}_\omega\ar@<2pt>[l] \ : \Tilde{m}^{k}}.$$
\end{theorem}

\begin{remark}
$f^k$ has a right adjoint $m^k$ and we have, for $F\in\Fo$, $m^k(F)=\Sigma^k m^{<1}(F)$.
\end{remark}

We conclude this subsection, by computing some useful $\Nil_k$-localisation. 

\begin{lemma}\label{nilkhst}
For $V\in\E$ and $M\in\U^{<k}$, $H^*(V)\otimes M$ is $nil_k$-closed. Furthermore, if $M$ is bounded in degree (with the bound possibly greater than $k$) the $nil_k$-localisation of $H^*(V)\otimes M$ is $H^*(V)\otimes M^{<k}$.
\end{lemma}

\begin{proof}
For $M$  bounded in degree, we have $$T_W(H^*(V)\otimes M)\cong T_W(H^*(V))\otimes T_W(M)\cong H^*(V)\otimes\F^{\Hom(W,V)}\otimes M,$$ and it's $H^*(W)$-comodule structure is the one of $T_W(H^*(V))$ tensor the identity of $T_W(M)\cong M$. Therefore, the projection from $H^*(V)\otimes M$ to $H^*(V)\otimes M^{<k}$ becomes an isomorphism when we apply $f^{<k}$. So, we only have to prove the first part of the Lemma.\\

For $M\in\U^{<k}$, let's prove that $H^*(V)\otimes M$ is $nil_k$-closed. $H^*(V)\otimes M$ admits an injective resolution $$0\rightarrow H^*(V)\otimes M\rightarrow I_0\rightarrow I_1\rightarrow ...$$
with $I_0$ and $I_1$ which are tensor products of $H^*(V)$ with some direct sum of $J(i)$ with $i<k$. Since, $f^{<k}$ is exact and $m^{<k}$ is left exact we have an exact sequence $$0\rightarrow l_k(H^*(V)\otimes M)\rightarrow l_k(I_0)\rightarrow l_k(I_1)\rightarrow...$$ But $I_0$ and $I_1$ are injectives and $nil_k$-reduced, they are therefore $nil_k$-closed. By the five lemma, we get that $H^*(V)\otimes M$ is also $nil_k$-closed.
\end{proof}

\subsection{$\Nil_k$-localisation of $T_V(K)$}

In this subsection, we define a shift functor on $\Fk{k}$ and prove that it is the ``$\Nil_k$-localisation of the $T$ functor''.

\begin{definition}\label{deltki}
For $V$ and $W$ in $\E$, let $\Delta_V\ :\ \Fk{k}\rightarrow\Fk{k}$ be the functor such that $\Delta_V F(W):=F(V\oplus W)$ and such that $\kappa_{\Delta_V(F),W}\ :\  \Delta_V(F)(W)\rightarrow\Delta_V(F)(W)\otimes_{\U^{<k}} H^{*<k}(W)$ is the composition of $\kappa_{F,V\oplus W}$ with $\id_{F(V\oplus W)}\otimes_{\U^{<k}} (\iota_W^*)^{<k}$, where $\iota_W$ is the canonical injection from $W$ to $V\oplus W$. 

\end{definition}

\begin{proposition}
The functor $\Delta_V$ is exact.
\end{proposition}

\begin{proposition}\label{bio}
Let $M$ be an unstable module and $V$ a finite dimensional vector space. There is a natural isomorphism $\Delta_V(f^{<k}(M))\cong f^{<k}(T_V(M)).$
\end{proposition}

\begin{proof}
We have $f^{<k}(T_V(M))(W)\cong\left[T_W(T_V(M))\right]^{<k}\cong T_{V\oplus W}(M)^{<k}\cong \Delta_V(f^{<k}(M))(W)$. Also, this isomorphism is compatible with the $H^{*<k}(W)$-comodule structures involved. This is because, by construction, the following diagram commutes:
$$\xymatrix{
T_{V\oplus W}(M)\ar[rr]^-{\kappa_{M,V\oplus W}}\ar[d]^{\cong} & &T_{V\oplus W}(M)\otimes H^*(V\oplus W)\ar[d]^-{id\otimes \iota_W^*}\\
T_W(T_V(M))\ar[rr]^-{\kappa_{T_V(M),W}}& &T_W(T_V(M))\otimes H^*(W).}$$ 
\end{proof}

\begin{lemma}\label{htM}
For $M\in\U^{<k}$, $T_V(H^*(W)\otimes M)$ is $nil_k$-closed.
\end{lemma}

\begin{proof}
Since the functor $T_V$ commutes with the tensor product, $$T_V(H^*(W)\otimes M)\cong T_V(H^*(W))\otimes T_V(M).$$ Then, since $M$ is bounded in degree $T_V(M)\cong M$, and therefore $T_V(H^*(W)\otimes M)\cong H^*(W)\otimes\F^{\Hom_{\E}(V,W)}\otimes M$. Thus, since $\F^{\Hom_{\E}(V,W)}\otimes M\in\U^{<k}$, Lemma \ref{nilkhst} implies that $T_V(H^*(W)\otimes M)$ is $nil_k$-closed.
\end{proof}

\begin{proposition}\label{tvnf}
If $M$ is $nil_k$-closed, $T_V(M)$ is also $nil_k$-closed.
\end{proposition}

\begin{proof}
Let $M$ be a $nil_k$-closed unstable module, and let $I_0$ and $I_1$ be direct sums in $\mathcal{U}$ of objects $H^*(W)\otimes J(i)$, such that $0\rightarrow M\rightarrow I_0\rightarrow I_1$ is exact. ($I_0$ and $I_1$ exist because of Corollary \ref{injreso}.) Since $f^{<k}$ is exact, $0\rightarrow f^{<k}(M)\rightarrow F_0\rightarrow F_1$ is exact for $F_0:=f^{<k}(I_0)$ and $F_1:=f^{<k}(I_1)$.\\

Since $m^{<k}$ is left exact, since $M$ is $nil_k$-closed and since, by Lemma \ref{nilkhst}, the $nil_k$-localisation of $H^*(W)\otimes J(i)$ is $H^*(W)\otimes J(i)^{<k}$, we have an exact sequence $0\rightarrow M\rightarrow I'_0\rightarrow I'_1$ where $I'_n$ is obtained from $I_n$ by replacing each $H^*(W)\otimes J(n)$ by $H^*(W)\otimes J(n)^{<k}$. Since $T_V$ is exact, the following is also exact: $$0\rightarrow T_V(M)\rightarrow T_V(I'_0)\rightarrow T_V(I'_1).$$

By applying the $nil_k$-localisation, we get the following commutative diagram:
$$\xymatrix{
0\ar[r]& T_V(M)\ar[r]\ar[d]^f& T_V(I'_0)\ar[r]\ar[d]^{m_0}& T_V(I'_1)\ar[d]^{m_1}\\
0\ar[r]&  l_k(T_V(M))\ar[r]& l_k(T_V(I'_0))\ar[r]& l_k(T_V(I'_1)),
}$$ where $m_0$, $m_1$ and $f$ are the $nil_k$-localisations. By Lemma \ref{htM}, $m_0$ and $m_1$ are isomorphisms. Hence, by the five lemma, $f$ is also an isomorphism which implies that $T_V(M)$ is $nil_k$-closed.
\end{proof}

\begin{corollary}\label{TVk}
Let $F\in\Fok{k}$, $m^{<k}(\Delta_ V(F))\cong T_V(m^{<k}(F))$.
\end{corollary}

\begin{proof}
We have $m^{<k}(\Delta_ V(F))\cong m^{<k}(\Delta_ V((f^{<k}\circ m^{<k})(F)))\cong (m^{<k}\circ f^{<k})(T_V(m^{<k}(F))$. The result follows from the fact that $T_V(m^{<k}(F))$ is $nil_k$-closed by Proposition \ref{tvnf}.
\end{proof}

\section{$\Nil_k$-localisation of objects from $\K$ and $K-\U$}

In this section, we prove that for $K\in\K$ and for $M\in K-\U$, $l_k(K)$ is an unstable algebra and $l_k(M)$ is an object in $l_k(K)-\U$. We start by proving that $m^{<k}$ is lax-monoidal.

\subsection{The functor $m^{<k}$ is lax-monoidal}

\begin{definition}
Let $\mu\ :\ F(*)\rightarrow F(*)\otimes F(*)$ be the morphism of graded unstable algebras defined by $\mu(x(k))=\bigoplus\limits_{i\leq k} x(i)\otimes x(i-k)$.\\

\end{definition}

From the coalgebra structure on $\A$, $F(*)\otimes F(*)$ gets a structure of left $\A$-module, and $\mu$ is then a morphism of left $\A$-modules.

\begin{proposition}\label{laxmon}
The functor $m^{<k}$ is lax-monoidal.
\end{proposition}

\begin{proof}
Let $F$ and $G$ be objects in $\Fk{k}$. From Proposition \ref{fufu}, we get that $m^{<k}(F)\cong\HomU(F(*),m^{<k}(F))$, which is isomorphic to $\Hom_{\Fk{k}}(f^{<k}(F(*)),F)$ by adjunction between $m^{<k}$ and $f^{<k}$. We define $\mu_{F,G}\ :\ m^{<k}(F)\otimes m^{<k}(G)\rightarrow m^{<k}(F\otimes G)$ in the following way. For $\alpha\in\Hom_{\Fk{k}}(f^{<k}(F(*)),F)$ and $\beta\in\Hom_{\Fk{k}}(f^{<k}(F(*)),G)$, $\mu_{F,G}(\alpha\otimes\beta)\in \Hom_{\Fk{k}}(f^{<k}(F(*)),F\otimes G)$ is the map which maps $x\in f^{<k}(F(*))(V)$ to $\bigoplus\limits_{(y,z)\in A} \alpha(y)\otimes\beta(z)$ for $A$ a part of $f^{<k}(F(*))(V)^2$ such that $f^{<k}(\mu)(x)=\bigoplus\limits_{(y,z)\in A} y\otimes z$. This map is a morphism of right $\A$-module, since $\mu$ is a morphism of left
$\A$-module.

The unitality, associativity and commutativity are straightforward.
\end{proof}

\subsection{Recollection about the functor $\Phi$}

In this subsection, we recall the definition of the functor $\Phi$ and recall it's connection to the definition of unstable algebras.

\begin{definition}\cite[I.7.2]{S1}
\begin{enumerate}
Let $\Phi$ be the functor from $\U$ to itself defined by $$(\Phi M)^n=\left\{\begin{array}{cc}
        0 & \text{for $n$ odd,}  \\
        M^{\frac{n}{2}} & \text{for $n$ even,} 
    \end{array}\right.$$ and
    $$\Sq^i\Phi x=\left\{\begin{array}{cc}
        0 & \text{for $i$ odd,}  \\
        \Phi(\Sq^{\frac{i}{2}}x) & \text{for $i$ even,} 
    \end{array}\right.$$ for all $x\in M$.
\end{enumerate}

\end{definition}

\begin{definition}
\begin{enumerate}
    \item For $M\in \U$, let $\lambda_M$ from $\Phi M$ to $M$ defined by $\Phi x\mapsto \Sq_0x$.
     \item For $K$ an algebra in $\U$, let $\sigma_K$ from $\Phi K$ to $K$ defined by $\Phi x\mapsto x^2$.

\end{enumerate}
\end{definition}

The following lemma follows directly from the definition of $\Phi$ and $\lambda_M$.

\begin{lemma}
For $M\in\U$, $\lambda_M$ is $\A$-linear.
\end{lemma}

We were interested in recalling the definition of $\Phi$ because we can characterise the fact that a $\A$-algebra is an unstable algebra in terms of $\lambda_K$ and $\sigma_K$. Indeed, an $\A$-algebra is an unstable algebra if and only if $\lambda_K=\sigma_K$.

\subsection{$\Nil_k$-localisation of $\K$ and $K-\U$ }

Since $f^{<k}$ is a strong monoidal functor and since $m^{<k}$ is lax-monoidal, $l_k$ is a lax-monoidal functor from $\U$ to itself. Therefore, for $K\in\K$, $l_k(K)$ gets a structure of $\A$-algebra. We want to show that this is a structure of unstable algebra.

\begin{lemma}\label{ouuf}
For $K\in\K$, the $\Nil_k$-localisation is a morphism of $\A$-algebras.
\end{lemma}

\begin{proof}

For $K\in\K$, the algebra structure of $K\cong\HomU(F(*),K)$ is the morphism which sends $\alpha\otimes\beta\in\HomU(F(*),K)^{\otimes 2}$, with $\alpha$ and $\beta$ of respective degrees $n$ and $k$, to the only element of degree $n+k$ of $\HomU(F(*),K)$ that send $x(n+k)$ to $\bigoplus\limits_{i\leq n+k} \alpha(x(i)).\beta(x(n+k-i))$. Indeed, since $\alpha(x(i))=0$ and $\beta(x(j))=0$ for $i\neq n$ and $j\neq k$, this is the morphism which sends $x(n+k)$ to $\alpha(x(n)).\beta(x(k))$. But $\alpha(x(n)).\beta(x(k))\in K$ is precisely the product of the elements $a$ and $b$ in $K$ represented by $\alpha$ and $\beta$ in $\HomU(F(*),K)$.\\

Therefore, we have the following commutative diagram: 

$$\xymatrix{\HomU(F(*),K)^{\otimes 2}\ar[r]\ar[d] & \Hom_{\Fk{k}}(f^{<k}(F(*)),f^{<k}(K))^{\otimes 2}\ar[d]\\
\HomU(F(*)\otimes F(*),K\otimes K)\ar[r]\ar[d] & \Hom_{\Fk{k}}(f^{<k}(F(*)\otimes F(*)),f^{<k}(K\otimes K))\ar[d]\\
\HomU(F(*),K)\ar[r] & \Hom_{\Fk{k}}(f^{<k}(F(*)),f^{<k}(K)),}$$
in which \begin{enumerate}
    \item the horizontal maps are induced by $f^{<k}$,
    \item the bottom vertical ones are induced by the multiplication of $K$ and the morphism $\mu$ from $F(*)$ to $F(*)^{\otimes 2}$,
    \item the first quadrant is commutative because $f^{<k}$ is a strong monoidal functor.
    
\end{enumerate}
Since the composition of the maps on the left is the multiplication of $K\cong\HomU(F(*),K)$ and the composition of the maps on the right is the multiplication of $l_k(K)\cong\Hom_{\Fk{k}}(f^{<k}(F(*)),f^{<k}(K))$, this proves the result.

\end{proof}

\begin{theorem}\label{princ1}
$l_k$ induces a functor from $\K$ to itself.
\end{theorem}

\begin{proof}

We already stated that $l_k(K)$ is an algebra in $\U$. Let us prove the unstability condition.\\ 

We have to prove that $\lambda_{l_k(K)}=\sigma_{l_k(K)}$. Since $K$ is an unstable algebra, we already know that $\lambda_K=\sigma_K$. We have the two following commutative diagrams: 
$$\xymatrix{\Phi K\ar[r]\ar[d]_{\lambda_K} & \Phi l_k(K)\ar[d]^{\lambda_{l_k(K)}} & \Phi K\ar[r]\ar[d]_{\sigma_K} & \Phi l_k(K)\ar[d]^{\sigma_{l_k(K)}}\\
K\ar[r] & l_k(K) & K\ar[r] & l_k(K),}$$
where the second one is commutative because of Lemma \ref{ouuf}. Applying the functor $f^{<k}$, we get two commutative diagrams whose horizontal maps are isomorphisms since we are applying $f^{<k}$ to some $\Nil_k$-localisation. Therefore, since $\lambda_K=\sigma_K$, we get that $f^{<k}(\lambda_{l_k(K)}-\sigma_{l_k(K)})=0$. This implies that for any $x\in\Ima(\sigma_{l_k(K)}-\lambda_{l_k(K)})$, the sub-module generated by $x$, $<x>$ is $k$-nilpotent. But $l_k(K)$ contains no non trivial $k$-nilpotent sub-modules so $x=0$, hence $\lambda_{l_k(K)}=\sigma_{l_k(K)}$.

\end{proof}

\begin{proposition}\label{princ2}
For $K\in\K$, $l_k$ induces a functor from $K-\U$ to $l_k(K)-\U$.
\end{proposition}

\begin{proof}
It is a direct consequence of the fact that $l_k$ is a lax-monoidal functor.
\end{proof}

\section{Central elements of an unstable algebra}

In this section, we recall the definition of central elements of an unstable algebra first introduced by Dwyer and Wilkerson in \cite{DW1}. We start by recalling the equivalence of categories between $\K/\Nil_1$ and some category of functors constructed in \cite{HLS2}, which will play an important role in the definition of central elements.

\subsection{$\Nil_1$-localisation of unstable algebras}

Since $\K$ is not abelian, one cannot define a localized category of $\K$ in the sense of \cite{G1}. In \cite{HLS1}, Henn, Lannes and Schwartz constructed a localized category $\K/\Nil_1$ with respect to the morphisms whose kernels and cokernels are in $\Nil_1$, in the sense of \cite{KS}. Then, the functor $f:=f^{<1}$ restricted to $\K$ factorises through a functor from $\K/\Nil_1$ to $\Fk{1}$. The authors of \cite{HLS1} identified the essential image of $f$ restricted to $\K$ and they deduced an equivalence of category between $\K/\Nil_1$ and some category of contravariant functors from $\E$ to the category of profinite sets. \\

\begin{definition}
A $2$-boolean algebra, is an algebra $B$ over $\F$, such that, for all $x\in B$, $x^2=x$.
\end{definition}

For any algebra $A\in\K$, $A^0$ is a $2$-boolean algebra. Therefore, for $K\in\K$, $T_V(K)$ is an unstable algebra, hence $T_V(K)^0$ is a $2$-boolean algebra. We can then use standard results on $2$-boolean algebras to study $f(K)$.\\

For $\mathcal{B}$ the category of $2$-boolean algebras and $B$ a $2$-boolean algebra, we consider $\Hom_{\mathcal{B}}(B,\F)$ the set of morphisms of $\F$-algebras from $B$ to $\F$. Since $B$ is the direct limit of its finite dimensional subalgebras $B_\alpha$, $\Hom_{\mathcal{B}}(B,\F)$ is the inverse limit of the $\Hom_{\mathcal{B}}(B_\alpha,\F)$ which are finite sets. $\Hom_{\mathcal{B}}(B,\F)$ inherits a structure of profinite set.

 \begin{proposition}\cite{HLS2}\label{booleprofin}
For $\mathcal{P}fin$ the category of profinite sets, the functor $\text{spec}\ :\ \mathcal{B}^{\text{op}}\rightarrow\mathcal{P}fin$, where $\text{spec}(B):=\text{Hom}_{\mathcal{B}}(B,\F)$, is an equivalence of categories whose inverse is the functor which sends $S$ to the algebra of continuous maps from $S$ to $\F$, $\F^S$.
\end{proposition}

In particular, for $K\in\K$, by adjunction $\Hom_{\mathcal{B}}(T_V(K)^0,\F)\cong\HomK(K,H^*(V))$ and this isomorphism is an isomorphism of profinite sets, where the structure of profinite set on $\HomK(K,H^*(V))$ comes from the fact that $K$ is the direct limit of the unstable sub-algebras of $K$ which are finitely-generated as $\A$-algebras. Then, $T_V(K)^0$ is isomorphic as a $2$-boolean algebra to $\F^{\HomK(K,H^*(V))}$. 

\begin{definition}
\begin{enumerate}
    \item Let $\PFVF$ be the category of functors from $(\E)^{\text{op}}$ to $\mathcal{P}fin$,
    \item let $\mathcal{L}$ be Lannes' linearization functor from $(\PFVF)^{\text{op}}$ to $\Fc$ defined by $\mathcal{L}(F)(V):=\F^{F(V)}$,
    \item let $g\ :\ \K^{\text{op}}\rightarrow\PFVF$ be the functor which sends $K$ to the functor $g(K)\ :\ V\mapsto\HomK(K,H^*(V))$.
\end{enumerate}
\end{definition}

We have a commutative diagram of functors:
$$\xymatrix{\K\ar[r]^-{g^{\text{op}}}\ar[d] & (\PFVF)^{\text{op}}\ar[d]^{\mathcal{L}}\\
\U \ar[r]^f &\Fc,}$$
where the functor from $\K$ to $\U$ is the forgetful functor. We denote by $\PFVF_\omega$ the full subcategory of $\PFVF$, whose objects are those whose image under $\mathcal{L}$ are in $\Fo$.

The functor $g$ has a unique factorisation of the following form:
$$\K^{\text{op}}\rightarrow (\K/\Nil)^{\text{op}}\rightarrow \PFVF_\omega\rightarrow \PFVF.$$

\begin{theorem}\cite[Theorem 1.5 of Part II]{HLS2}
The functor from $\K/\Nil$ to $\PFVF_\omega$ induced by $g^{\text{op}}$ is an equivalence of categories.
\end{theorem}

The following lemma will be of importance in the following.

\begin{lemma}\label{coco limite}
The functor $g$ turns injections in $\K$ into surjections in $\PFVF$ and finite inverse limits in $\K$ into direct limits in $\PFVF$.
\end{lemma}

\begin{proof}
Since $f$ is exact, $f$ sends injections into injections and commutes with finite inverse limits. The result is then a consequence of the isomorphism $f(K)\cong \F^{g(K)}$.
\end{proof} 
\color{black}

\begin{definition}
For $K\in\K$, let $\mathcal{S}(K)$ be the category whose objects are pairs $(V,\phi)$ with $V\in\E$ and $\phi\in\HomK(K,H^*(V))$, and whose morphisms between $(V,\phi)$ and $(W,\psi)$ are linear maps $\alpha$ from $V$ to $W$ such that $\alpha^*\psi=\phi$.\\
Let also $\Skop{K}$ denote the category of contravariant functor from $\Skop{K}$ to $\E$.
\end{definition}

\subsection{Connected components of $T_V(K)$}

For $K$ an unstable algebra and $M$ a $K$-module in $\U$, we recall the definition of the connected components of $T_V(K)$ and $T_V(M)$ over $T_V(K)^0$, which is reviewed in \cite{heard2019depth} and \cite{heard2020topological}. Such a decomposition exists for any graded module or algebra over a $2$-boolean algebra.

\begin{lemma}
For $K$ an unstable algebra which is finitely-generated as an algebra over $\A$,\newline $\text{Hom}_\mathcal{K}(K,H^*(V))$ is finite. 
\end{lemma}

Let us first notice that, for $K\in\K$ and $M\in K-\U$ the category of $K$-modules in $\U$, $T_V(K)$ is in $\K$ and $T_V(M)$ in $T_V(K)-\U$. This is because the functor $T_V$ commutes with tensor products (Proposition \ref{muche}). In particular, $T_V(K)$ is a $T_V(K)^0$-algebra and $T_V(M)$ is a $T_V(K)^0$-module.

\begin{definition}
\begin{enumerate}
    \item For $V\in\E$ and $\phi \in \text{Hom}_\mathcal{K}(K,H^*(V))$, let $$T_{(V,\phi)}(K):=T_V(K)\otimes_{T_V(K)^0}\F(\phi),$$ where the structure of $T_V(K)^0$-module on $\F(\phi)$ is induced by the morphism from $T_V(K)^0$ to $\F$ adjoint to $\phi$.
    \item Let also $T_{(V,\phi)}(M):=T_V(M)\otimes_{T_V(K)^0}\F(\phi)$.
\end{enumerate}

\end{definition}

\begin{notation}
We denote respectively by $\zeta_{K,(V,\phi)}$ and $\zeta_{M,(V,\phi)}$ the canonical projections from $T_V(K)$ and $T_V(M)$ to $T_{(V,\phi)}(K)$ and $T_{(V,\phi)}(M)$. 
\end{notation}

\begin{proposition}\cite[Equation (2.6)]{heard2019depth}\label{prumi}
\begin{enumerate}
    \item For $K$ an unstable algebra and $V\in\E$, with $\HomK(K,H^*(V))$ finite, we have the following natural isomorphism of unstable algebras $$T_V(K)\cong \prod\limits_{\phi\in \text{Hom}_\mathcal{K}(K,H^*(V))}T_{(V,\phi)}(K).$$
    \item Also, for $M\in K-\U$, we have the following natural isomorphism of unstable modules $$T_V(M)\cong \bigoplus\limits_{\phi\in \text{Hom}_\mathcal{K}(K,H^*(V))}T_{(V,\phi)}(M).$$
\end{enumerate}

\end{proposition}

\poubelle{
In the following, it will be useful to have a canonical identification of $T_{(V,\phi)}(K)$ and $T_{(V,\phi)}(M)$ as sub-algebras of $T_V(K)$ and $T_V(M)$. In the following, for $\phi\in\HomK(K,H^*(V)$, we denote by $\delta_\phi\in T_V(K)^0\cong\F^{\HomK(K,H^*(V))}$ the characteristic function of $\phi$ and by $\delta_\phi\times T_V(K)$ (respectively $\delta_\phi\times T_V(M)$) the image of the multiplication by $\delta_\phi$ in $T_V(K)$ (respectively $T_V(M)$).

\begin{lemma}\label{technity}
For $K$ an unstable algebra, $M\in K-\U$, $V\in\E$ and $\phi\in\HomK(K,H^*(V))$, the isomorphisms of Proposition \ref{prumi} identify $T_{(V,\phi)}(K)$ and $T_{(V,\phi)}(M)$ with $\delta_\phi \times T_V(K)$ and $\delta_\phi \times T_V(M)$. Furthermore, $\zeta_{K,(V,\phi)}$ and $\zeta_{M,(V,\phi)}$ identify with the multiplication by $\delta_\phi$.
\end{lemma}

\begin{proof}

We give the proof for $K$; the one for $M$ is similar.\\

We start by noticing that, for any $x\in T_V(K)$, $\zeta_{K,(V,\phi)}(\delta_\phi\times x)=\zeta_{K,(V,\phi)}(\delta_\phi)\times\zeta_{K,(V,\phi)}(x)$. Since $\zeta_{K,(V,\phi)}(\delta_\phi)=1$, this is equal to $\zeta_{K,(V,\phi)}(x)$. Hence, $\zeta_{K,(V,\phi)}$ factorises through its own restriction to $\delta_\phi\times T_V(K)$. The result follows from the fact that $\zeta_{K,(V,\phi)}(x)=0$ if and only if $\delta_\phi\times x=0$.
\end{proof}}

\begin{lemma}\label{upussu}
For $K$ an unstable algebra, $V\in\E$ such that $\HomK(K,H^*(V))$ is finite, $M\in K-\U$, and $\phi\in\HomK(K,H^*(V))$,\begin{enumerate}
    \item if $K$ is $nil_k$-closed, $T_{(V,\phi)}(K)$ is $nil_k$-closed,
    \item if $M$ is $nil_k$-closed, $T_{(V,\phi)}(M)$ is $nil_k$-closed.
\end{enumerate} 
\end{lemma}

\begin{proof}
By Proposition \ref{prumi}, $T_{(V,\phi)}(K)$ is a direct summand of $T_V(K)$. Also, by Proposition \ref{tvnf} $T_V(K)$ is $nil_k$-closed. Since $l_k$ is additive, this implies that $T_{(V,\phi)}(K)$ is $nil_k$-closed. The proof for $M$ is similar.
\end{proof}

\subsection{Injective objects in $K-\U$}

\begin{definition}\cite[]{BSMF_1995__123_2_189_0}
For $\phi\in\HomK(K,H^*(V))$. Let $H^*(V)^\phi\otimes-$ be the functor from $T_{(V,\phi)}(K)-\U$ to $K-\U$ which maps $M\in T_{(V,\phi)}(K)-\U$ to $H^*(V)\otimes M$ with the $K$-module structure induced by its $H^*(V)\otimes T_{(V,\phi)}(K)$-module structure and $\eta_{K,(V,\phi)}\ :\ K\rightarrow H^*(V)\otimes T_{(V,\phi)}(K)$. 
\end{definition}

\begin{proposition}\cite[]{BSMF_1995__123_2_189_0}\label{adjK}
The functor $T_{(V,\phi)}$ is left adjoint to $H^*(V)^\phi\otimes-$.

\end{proposition}

In the article \cite{BSMF_1995__123_2_189_0}, the authors exhibit injective cogenerators of the category $K-\U$.

\begin{definition}\cite[]{BSMF_1995__123_2_189_0}\label{injK}
Let $J_K(n)$ be the object of $K-\U$ determined up to isomorphism by the natural isomorphism in $M$, $\Hom_{K-\U}(M,J_K(n))\cong (M^n)^\sharp$.
\end{definition}

\begin{theorem}
$(J_K(n))_{n\in\N}$ is a family of injective cogenerators of $K-\U$.
\end{theorem}

\begin{definition}\cite[]{Henn+1996+189+215}
Let $I_{(V,\phi)}(n):=H^*(V)\otimes^\phi J_{T_{(V,\phi)(K)}(n)}$.
\end{definition}

\begin{proposition}\cite[]{Henn+1996+189+215}\label{injmodK}
For $n\in\N$ and $\phi\in\HomK(K,H^*(V))$, $I_{(V,\phi)}(n)$ is injective and we have the following natural isomorphism in $M\in K-\U$, $\Hom_{K-\U}(M,I_{(V,\phi)}(n))\cong (T_{(V,\phi)}(M)^n)^\sharp$.
\end{proposition}

\begin{proof}
This is a direct consequence of Definition \ref{injK} and Proposition \ref{adjK}.
\end{proof}

\begin{proposition}\cite[1.8]{Henn+1996+189+215}
For $K$ a noetherian unstable algebra, $\phi\in\HomK(K,H^*(V))$ and $n\in\N$, $I_{(V,\phi)}(n)$ is finitely-generated as a $K$-module.
\end{proposition}

\begin{corollary}\label{gazougi}
For $K$ a noetherian unstable algebra and $M\in K-\U$ finitely-generated as a $K$-module, $\Hom_{K-\U}(M,I_{(V,\phi)}(n))$ is finite, therefore $$T_{(V,\phi)}(M)^n\cong \Hom_{K-\U}(M,I_{(V,\phi)}(n))^\sharp.$$
\end{corollary}

\begin{example}
For $K\in\K$, $M\in K-\U$, $V\in\E$ and $\phi\in\HomK(K,H^*(V))$, $T_{(V,\phi)}(M)^0\cong \Hom_{K-\U}(M,H^*(V,\phi))^\sharp$ where $H^*(V,\phi)$ is the algebra $H^*(V)$ with the $K$-algebra structure induced by $\phi$.
\end{example}

\begin{definition}
For $M\in K-\U$ finitely-generated, let $\Hom_{K-\U}(M,I_{(\_,\_)}(n))$ be the contravariant functor from $\mathcal{S}(K)$ to $\E$ which maps $(V,\phi)$ to $\Hom_{K-\U}(M,I_{(V,\phi)}(n))$.
\end{definition}

\begin{notation}
For $n=0$, since $I_{(V,\phi)}(0)=H^*(V,\phi)$ we will denote $\Hom_{K-\U}(M,I_{(\_,\_)}(0))$ by $\Hom_{K-\U}(M,H^*(\_,\_))$.
\end{notation}

\begin{proposition}\label{keykey}
For $K$ a noetherian unstable algebra and $M$ finitely-generated as a $K$-module, $f^k(M)(V)$ is naturally isomorphic to $\bigoplus\limits_{\phi\in\HomK(K,H^*(V))}\Hom_{K-\U}(M,I_{(V,\phi)}(k))^{\sharp}$.
\end{proposition}

\begin{proof}
This is a direct consequence of Corollary \ref{gazougi} and Proposition \ref{prumi}.
\end{proof}

\subsection{Definition and first properties of central elements}

The notion of a central element of an unstable algebra $K$ was defined by Dwyer and Wilkerson in \cite{DW1}; they used it in \cite{DW12} to exhibit the only exotic finite loop space at the prime $2$. The centre of $K$ has been studied in detail in \cite{heard2019depth} and \cite{heard2020topological}.

The aim of this subsection is to recall some known facts about central elements of an unstable algebra.\\
 
 \begin{notation}
 Let $N$ be an unstable module, $V$ a finite dimensional vector space, $K$ an unstable algebra, $\phi\in\HomK(K,H^*(V))$ and $M\in K-\U$. Denote by:
\begin{itemize}  
    \item  $\eta_{N,V}\ :\ N\rightarrow T_V(N)\otimes H^*(V)$ the unit of the adjunction between $T_V$ and $-\otimes H^*(V)$;
    \item  $\rho_{N,V}$ the following composition $$\rho_{N,V}\ :\ N\overset{\eta_{N,V}}{\longrightarrow} T_V(N)\otimes H^*(V)\overset{id\otimes \epsilon_{V}}{\longrightarrow} T_V(N), $$
where $\epsilon_V$ denote the augmentation of $H^*(V)$;
     \item  $\rho_{K,(V,\phi)}$ the composition of $\rho_{K,V}$ with the projection onto $T_{(V,\phi)}(K)$;
     \item  $\rho_{M,(V,\phi)}$ the composition of $\rho_{M,V}$ with the projection onto $T_{(V,\phi)}(M)$.
\end{itemize}
 \end{notation} 

\begin{remark}
The morphism $\rho_{N,V}$ identifies with $T_{\iota_0^V}(N)\ :\ N\cong T_0(N)\rightarrow T_V(N)$, the morphism induced by naturality of $T_V(N)$ with respect to $V$ by $\iota_0^V$, the injection from $0$ to $V$.
\end{remark}

\begin{definition}\label{centri}
\begin{enumerate}
    \item For $K$ an unstable algebra and $\phi\in\HomK(K,H^*(V))$, the pair $(V,\phi)$ is said to be central if $\rho_{K,(V,\phi)}\ :\ K\rightarrow T_{(V,\phi)}(K)$ is an isomorphism. 
    \item For $M\in K-\U$, $(V,\phi)$ is said to be $M$-central if $\rho_{M,(V,\phi)}\ :\ M\rightarrow T_{(V,\phi)}(M)$ is an isomorphism. 
    \item Let $\textbf{C}(K)$ be the set of central elements of $K$ and let $\textbf{C}_K(M)$ be the set of $M$-central elements.
\end{enumerate}

\end{definition}

The classical example and first motivation for studying the centre of an unstable algebra is the example of $H^*(G)$, the cohomology of a group $G$. The details can be found in \cite{Henn2001}.

\begin{example}
For $G$ a discrete group or a compact Lie group, $\text{Hom}_\mathcal{K}(H^*(G),H^*(V))\cong\F\left[\text{Rep}(V,G)\right]$, where $\text{Rep}(V,G)$ denote the conjugacy classes of morphisms from $V$ to $G$.\\

Let $\rho$ represent a conjugacy class in $\text{Rep}(V,G)$. We consider the morphism $V\times C_G(\rho)\rightarrow G$, where $C_G(\rho)$ denote the centraliser in $G$ of the image of $\rho$, which sends $(v,g)$ to $\rho(v)\cdot g$. It induces a morphism from $H^*(G)\rightarrow H^*(V)\otimes H^*(C_G(\rho))$. By adjunction, it gives us a morphism $T_V(H^*(G))\rightarrow H^*(C_G(\rho))$ which depends only on the conjugacy class of $\rho$. This morphism induces an isomorphism between $T_{(V,\rho)}(H^*(G))$ and $H^*(C_G(\rho))$, and $$\rho_{H^*(G),(V,\rho)}\ :\ H^*(G)\rightarrow T_{(V,\rho)}(H^*(G))\cong H^*(C_G(\rho))$$ is the morphism induced by the injection $C_G(\rho)\hookrightarrow G$.\\

Hence, $(V,\rho)$ is central if and only if the injection $C_G(\rho)\hookrightarrow G$ induces an isomorphism in cohomology.
\end{example}

\begin{definition}
Let $K$ be an unstable algebra, $K$ is connected if the unit $\eta\ :\ \F\rightarrow K$ is an isomorphism. If $K$ is connected, $K$ admit a unique augmentation $\epsilon_K\ :\ K\rightarrow\F$.
\end{definition}

\begin{remark}

The components $T_{(V,\phi)}(K)$ are connected by construction. Hence, if $K$ is not connected, $\textbf{C}(K)=\emptyset$.

\end{remark}

\begin{example}
The functor $T_0$ is the identity, hence, if $K$ is connected, $T_{(0,\epsilon_K)}(K)\cong K$, for $\epsilon_K\ :\ K\rightarrow \F$ the unique non-trivial morphism of unstable algebras from $K$ to $\F$. Hence $(0,\epsilon_K)$ is central.
\end{example}

\begin{notation}
Let $\epsilon_{K,V}$, be the composition of $\epsilon_K$ with the injection from $\F$ to $H^*(V)$.
\end{notation}

For $K$ a connected unstable algebra, $I(K)$ denotes the augmentation ideal of $K$; the module of indecomposable elements of $K$ is the unstable module defined by $\mathcal{Q}(K):=I(K)/I(K)^2$. An unstable module $M$ is said to be locally finite if, for all $x\in M$, $\A x$ is finite.\\

In \cite{AST_1990__191__97_0}, Dwyer and Wilkerson give several criterion to determine central elements of an unstable algebra $K$, in the case where $\mathcal{Q}(K)$ is locally finite. 

\begin{proposition} \cite[Proof of Theorem 3.2]{AST_1990__191__97_0}\label{indecompi}

Let $K$ be a connected unstable algebra such that $\mathcal{Q}(K)$ is locally finite as an unstable module, then  $(V,\epsilon_{K,V})$ is central for all $V\in\E$.

\end{proposition}

In particular, if $K$ is a connected, noetherian, unstable algebra, then $(V,\epsilon_{K,V})$ is central for all vector spaces $V$.

\begin{proposition} \cite[Proposition 3.4]{DW1}\label{DWC}

Let $K$ be a connected unstable algebra such that $\mathcal{Q}(K)$ is locally finite. Then, for\newline $\phi\in\HomK(K,H^*(V))$, $(V,\phi)$ is central if and only if there exists a morphism from $K$ to $K\otimes H^*(V)$ such that the following diagram commutes:
$$\xymatrix{ & & K\\
K\ar[rru]^{\text{id}}\ar[rr]\ar[rrd]_\phi & & K\otimes H^*(V)\ar[u]_-{\id\otimes\epsilon_{H^*(V)}}\ar[d]^-{\epsilon_K\otimes\id} & \\
& & H^*(V).}$$

\end{proposition}

\begin{corollary}\cite{DW1}\label{DWCT}
Let $K$ be a connected unstable algebra such that $\mathcal{Q}(K)$ is locally finite. For $\phi\in\HomK(K,H^*(V))$, $(V,\phi)$ is central if and only if $K$ has a structure of $H^*(V)$-comodule $\kappa$ in $\K$, such that the following diagram commutes:
$$\xymatrix{ 
K\ar[rr]^\kappa\ar[rrd]_\phi & & K\otimes H^*(V)\ar[d]^-{\epsilon_K\otimes\id} & \\
& & H^*(V).}$$
\end{corollary}

In particular, this implies:

\begin{corollary}\label{incl}
Let $K$ be an unstable algebra such that $\mathcal{Q}(K)$ is locally finite, then for  $\phi\in\textbf{C}(K)$ and $\alpha\ :\ V\rightarrow E$ a morphism in $\E$, $(V,\alpha^*\circ \phi)\in\textbf{C}(K)$.
\end{corollary}

\begin{proposition}
For any morphism of unstable algebra $\phi\ :\ H^*(W)\rightarrow H^*(V)$, $(V,\phi)\in\textbf{C}(H^*(W))$.
\end{proposition}

\begin{proof}
We consider $\nabla_W^*$, the coalgebra structure on $H^*(W)$ induced by $\nabla_W\ :\ W\oplus W \rightarrow W$. Then, the composition $(\id\otimes\phi)\circ \nabla^*_W$ is a $H^*(V)$-comodule structure satisfying the hypothesis of Corollary \ref{DWCT}. Therefore $(V,\phi)$ is central.
\end{proof}

\subsection{Centrality for $K$-algebras}

Let $K'$ be an object in $K\downarrow\K$, in this subsection we study $K'$-central elements in $\mathcal{S}(K)$. This will give us some first examples of $M$-centrality, for $M\in K-\U$, the case of algebras being a lot easier than the case where $M$ is just an object in $K-\U$.

For $\gamma\ :\ K\rightarrow K'$ a morphism in $\K$, $K'$ inherits a structure of unstable $K$-module. Let's assume that $K$ and $K'$ are both finitely-generated as algebras over $\A$. From Proposition \ref{prumi}, $T_V(K')$ admit two decompositions $T_V(K')\cong\prod\limits_{\psi\in\HomK(K,H^*(V))}T_{(V,\psi)}(K')$ and $T_V(K')\cong\prod\limits_{\phi\in\HomK(K',H^*(V))}T_{(V,\phi)}(K')$. Both are isomorphisms of unstable $T_V(K)$-algebras. We get an isomorphism of unstable $T_V(K)$-algebras $$\prod\limits_{\psi\in\HomK(K,H^*(V))}T_{(V,\psi)}(K')\cong\prod\limits_{\phi\in\HomK(K',H^*(V))}T_{(V,\phi)}(K').$$

\begin{lemma}
Under this isomorphism, $T_{(V,\psi)}(K')\cong\prod\limits_{\overset{\phi\in\HomK(K',H^*(V))\ ;}{\phi\circ\gamma=\psi}}T_{(V,\phi)}(K')$. 
\end{lemma}

\begin{proposition}\label{youpcentr}
For $K'$ a connected unstable algebra, finitely-generated as an algebra over $\A$, for $\gamma\ :\ K\rightarrow K'$ a morphism in $\K$ with $K$ also finitely-generated as a $\A$-algebra and for $\psi\in\HomK(K,H^*(V))$, $(V,\psi)$ is $K'$-central if and only if $\psi$ has a unique element $\phi$ in it's preimage under $\gamma^*$ in $\HomK(K',H^*(V))$ and if $(V,\phi)\in\textbf{C}(K')$. 
\end{proposition}

\begin{proof}
$K'$ being connected, the morphism $$\rho_{K',(V,\psi)}\ :\ K'\rightarrow T_{(V,\psi)}(K')\cong\prod\limits_{\overset{\phi\in\HomK(K',H^*(V))\ ;}{\gamma^*\phi=\psi}}T_{(V,\phi)}(K')$$ may be an isomorphism, only if $T_{(V,\psi)}(K')\cong\prod\limits_{\overset{\phi\in\HomK(K',H^*(V))\ ;}{\gamma^*\phi=\psi}}T_{(V,\phi)}(K')$ is connected. This is the case if and only if $\{\phi\in\HomK(K',H^*(V))\ ;\ \gamma^*\phi=\psi\}$ is a singleton. In this case, by definition, $\rho_{K',(V,\psi)}=\rho_{K',(V,\phi)}$ is an isomorphism if and only if $(V,\phi)$ is central.
\end{proof}

As a trivial special case, we get that, for $K\in\K$ and $(V,\psi)$ in $\mathcal{S}(K)$, $(V,\psi)$ is central if and only if it is $K$-central.

\begin{corollary}
For $K\in\K$, $K$ is an object in $K-\U$ and as such, $\textbf{C}(K)=\textbf{C}(K;K)$.
\end{corollary}

\begin{example}\label{gluglu}
Let $K$ be an algebra, finitely-generated as an $\A$-algebra, and let $\gamma\ :\ K\rightarrow H^*(W)$. Then, $\phi\in\HomK(K,H^*(V))$ is $H^*(W)$-central if and only if $\phi$ has a unique element in it's preimage under $\gamma^*$ in $\HomF(V,W)$. In particular, if $K=H^*(W)^G$ for $G$ a subgroup of $Gl(W)$, $\left[\phi\right]\in\HomF(V,W)/G$ is $H^*(W)$-central if and only if $\left[\phi\right]$ is a singleton. This is the case if and only if for all $x\in\Ima(\phi)$ and for all $g\in G$, $gx=x$. 
\end{example}

\section{Central elements and the nilpotent filtration}

In this section we are interested in the relationship between $\textbf{C}(K)$ and the nilpotent filtration, for $K\in \K$ noetherian.\\

More precisely, for $\phi\in\HomK(K,H^*(W))$, we say that $(W,\phi)$ is central away from $\Nil_k$ if $f^{<k}(\rho_{K,(W,\phi)})\ : f^{<k}(K)\rightarrow f^{<k}(T_{(W,\phi)}(K))$ is an isomorphism. If $(W,\phi)$ is central, $\rho_{K,(W,\phi)}$ is an isomorphism, hence it is central away from $\Nil_k$ for all $k$.\\ 

Since $f^{<k}(K)(W)$ identifies with $(f^{<k+1}(K)(W))^{<k}$, if $(W,\phi)$ is central away from $\Nil_{k+1}$ it is also central away from $\Nil_k$.\\ 

In this section, we will give a strategy to determine $\textbf{C}(K)$ by studying obstructions, for a pair $(W,\phi)$ central away from $\Nil_k$, to be also central away from $\Nil_{k+1}$.

\subsection{Centrality away from $\Nil_k$}

\begin{definition}
For $k\in\N^*$, $K\in\K$ and $M\in K-\U$, \begin{enumerate}
    \item let $\textbf{C}_k(K)$ be the following set: $$\textbf{C}_k(K):=\{(W,\phi)\text{ with } W\in\E\text{ and } \phi\in\HomK(K,H^*(W))\ |\ f^{<k}(\rho_{K,(W,\phi)}) \text{ is an isomorphism}\},$$
    \item and let $\textbf{C}_k(M;K)$ be: $$\textbf{C}_k(M;K):=\{(W,\phi)\text{ with } W\in\E\text{ and } \phi\in\HomK(K,H^*(W))\ |\ f^{<k}(\rho_{M,(W,\phi)}) \text{ is an isomorphism}\}.$$
\end{enumerate}
\end{definition}

\begin{remark}
For $V$ in $\E$, as a morphism of $\F$-vector spaces, the natural transformation $f^{<k}(\rho_{K,(W,\phi)})_V$ and $f^{<k}(\rho_{M,(W,\phi)})_V$ are the restrictions to degrees lower than $k$ of the maps $T_V(\rho_{K,(W,\phi)})$ and $T_V(\rho_{M,(W,\phi)}$. Therefore, we have the following towers of inclusions: $$\textbf{C}(K)\subset...\subset\textbf{C}_k(K)\subset ... \subset\textbf{C}_2(K)\subset\textbf{C}_1(K)$$
and $$\textbf{C}(M;K)\subset...\subset\textbf{C}_k(M;K)\subset ... \subset\textbf{C}_2(M;K)\subset\textbf{C}_1(M;K).$$ 
\end{remark}

\begin{proposition}\label{princ3}
Let $K\in\K$ and $M\in K-\U$, then $$\textbf{C}(K)=\bigcap\limits_{k\in\N^*}\textbf{C}_k(K)$$
and $$\textbf{C}(M;K)=\bigcap\limits_{k\in\N^*}\textbf{C}_k(M;K).$$
\end{proposition}

\begin{proof}
Since the functor $T_V$ is exact, and since $T_V(K)=0$ if and only if $K=0$, $T_V(\rho_{K,(W,\phi)})$ is an isomorphism if and only if $\rho_{K,(W,\phi)}$ is. Furthermore, if $f^{<k}(\rho_{K,(W,\phi)})_V=T_V(\rho_{K,(W,\phi)})^{<k}$ is an isomorphism for all $k\in\N$, $T_V(\rho_{K,(W,\phi)})$ is an isomorphism in all degrees therefore it is an isomorphism. The proof for $M$ is similar.
\end{proof}

For $K\in\K$, let us give an alternative description of $\textbf{C}_k(K)$. Any morphism $\phi$ from $K$ to $H^*(W)$ factorises uniquely as $K\rightarrow l_k(K)\overset{l_k(\phi)}{\rightarrow} H^*(W)$, where the first morphism is the localisation away from $\Nil_k$. We get from this the following proposition.

\begin{proposition}\label{nilfer}
For $K\in\K$ and $(W,\phi)\in\mathcal{S}(K)$, $(W,\phi)\in\textbf{C}_k(K)$, if and only if $(W,l_k(\phi))\in\textbf{C}(l_k(K))$.
\end{proposition}

\begin{proof}
If $\phi\in \textbf{C}_k(K)$, then $f^{<k}(\rho_{K,(W,\phi)})\ :\ f^{<k}(K)\rightarrow f^{<k}(T_{(W,\phi)}K)$ is an isomorphism. Then, applying the functor $m^{<k}$ to $f^{<k}(\rho_{K,(W,\phi)})\ :\ f^{<k}(K)\overset{\cong}{\rightarrow} f^{<k}(T_{(W,\phi)}K)$, we get that $\rho_{l_k(K),(W,l_k(\phi))}\ :\ l_k(K)\rightarrow T_{(W,l_k(\phi))}(l_k(K))$ is an isomorphism in $\U/\Nil_k$. Since $l_k(K)$ is $nil_k$-closed by definition, and since $T_{(W,l_k(\phi))}(l_k(K))$ is $nil_k$-closed by Lemma \ref{upussu}, this is an isomorphism in $\K$. Thus, $(W,l_k(\phi))$ is central.\\

Conversely, we have the following commutative diagram:
$$\xymatrix{f^{<k}(K)\ar[r] \ar[d] & f^{<k}(T_{(W,\phi)}(K))\ar[d]\\
f^{<k}(l_k(K))\ar[r]  & f^{<k}(T_{(W,l_k(\phi))}(l_k(K))),}$$
Where the vertical maps are isomorphisms induced by the localisation away from $\Nil_k$, and the horizontal ones are induced $\rho_{K,(W,\phi)}$ and $\rho_{l_k(K),(W,l_k(\phi))}$. Therefore, if $l_k(\phi)$ is central, $f^{<k}(\rho_{K,(W,\phi)})$ is an isomorphism so that, $\phi\in\textbf{C}_k(K)$.
\end{proof}

Similarly, we have the following.

\begin{proposition}\label{nilfer2}
For $K\in\K$, $M\in K-\U$ and $(W,\phi)\in\mathcal{S}(K)$, $(W,\phi)\in\textbf{C}_k(M;K)$, if and only if $(W,l_k(\phi))\in\textbf{C}(l_k(M);l_k(K))$.
\end{proposition}

We will now describe the obstruction for an element $(W,\phi)\in\textbf{C}_k(K)$ to be in $\textbf{C}_{k+1}(K)$.\\

Let $K$ be an unstable algebra and let $\lambda_k\ :\ l_{k+1}(K)\rightarrow l_k(K)$ be the $nil_{k}$-localisation of $l_{k+1}(K)$. We have an exact sequence in $K-\U$ 
$$0\rightarrow\ker(\lambda_k)\rightarrow l_{k+1}(K)\overset{\lambda_k}{\rightarrow}l_k(K)\rightarrow \coker(\lambda_k)\rightarrow 0.$$ We notice that $\ker(\lambda_k)$ is an ideal of $l_{k+1}(K)$, hence $\ker(\lambda_k)$ is an object in $K-\U$. Also $l_{k+1}(K)$ and $l_k(K)$ are $K$-algebras and therefore $K$-modules. Finally, $\lambda_k$ is a morphism of $K$-modules (this is because the localisation away from $\Nil_k$ of $K$ is canonically isomorphic to the localisation away from $\Nil_k$ of its localisation away from $\Nil_{k+1}$) therefore $\coker(\lambda_k)$ is also in $K-\U$. By naturality of $\rho_{M,(W,\phi)}$ with respect to $M\in K-\U$ (and using that $K$ is an object in $K-\U$), we get the following commutative diagram in $K-\U$, with exact rows:

$$
\resizebox{\displaywidth}{!}{%
\xymatrix{0\ar[r]  &\ker(\lambda_k)\ar[r] \ar[d]^-{\rho_{\ker(\lambda_k),(W,\phi)}} & l_{k+1}(K)\ar[r] \ar[d]^-{\rho_{l_{k+1}(K),(W,\phi)}} & l_k(K)\ar[r] \ar[d]^-{\rho_{l_k(K),(W,\phi)}} & \coker(\lambda_k)\ar[r]\ar[d]^-{\rho_{\coker(\lambda_k),(W,\phi)}} & 0\\
0\ar[r]  & T_{(W,\phi)}(\ker(\lambda_k))\ar[r]  & T_{(W,\phi)}(l_{k+1}(K))\ar[r]  & T_{(W,\phi)}(l_k(K))\ar[r] & T_{(W,\phi)}(\coker(\lambda_k))\ar[r] & 0.}}$$

By applying the exact functor $f^{<k+1}$, we get the following commutative diagram in $K-\U$ with exact rows.

$$\resizebox{\displaywidth}{!}{%
\xymatrix@C=1em{0\ar[r]  & f^{<k+1}(\ker(\lambda_k))\ar[r] \ar[d]^-{f^{<k+1}(\rho_{\ker(\lambda_k),(W,\phi)})} & f^{<k+1}(K)\ar[r] \ar[d]^-{f^{<k+1}(\rho_{K,(W,\phi)})} & f^{<k+1}(l_k(K))\ar[r] \ar[d]^-{f^{<k+1}(\rho_{l_k(K),(W,\phi)})} & f^{<k+1}(\coker(\lambda_k))\ar[r]\ar[d]^-{f^{<k+1}(\rho_{\coker(\lambda_k),(W,\phi)})} & 0\\
0\ar[r]  & f^{<k+1}(T_{(W,\phi)}(\ker(\lambda_k)))\ar[r]  & f^{<k+1}(T_{(W,\phi)}(l_{k+1}(K)))\ar[r]  &f^{<k+1}(T_{(W,\phi)}(l_k(K)))\ar[r] & f^{<k+1}(T_{(W,\phi)}(\coker(\lambda_k)))\ar[r] & 0.}}$$

\begin{theorem}\label{gaga}
Let $K$ be an unstable algebra and let $\phi\in\HomK(K,H^*(W))$. Then, the two following conditions are equivalent:
\begin{enumerate}
    \item $(W,\phi)\in\textbf{C}_{k+1}(K)$,
    \item $(W,\phi)\in\textbf{C}_k(K)$, $(W,\phi)\in\textbf{C}_{k+1}(\ker(\lambda_k);K)$ and  $(W,\phi)\in\textbf{C}_{k+1}(\coker(\lambda_k);K)$.
\end{enumerate}
\end{theorem}

\begin{proof}
We assume that $(W,\phi)\in\textbf{C}_{k+1}(K)$ then, $f^{<k+1}(\rho_{K,(W,\phi)})$ is an isomorphism. By restricting to degrees less than $k$, we get that $f^{<k}(\rho_{K,(W,\phi)})$ also is an isomorphism. Therefore, by Proposition \ref{nilfer}, $f^{<k+1}(\rho_{l_k(K),(W,\phi)})=f^{<k+1}(\rho_{l_k(K),(W,l_k(\phi))})$ is also an isomorphism. By the five lemma, we conclude that $f^{<k+1}(\rho_{\ker(\lambda_k),(W,\phi)})$ and $f^{<k+1}(\rho_{\coker(\lambda_k),(W,\phi)})$ are also isomorphisms.\\

Conversely, if $(W,\phi)$ is in $\textbf{C}_k(K)$ and if $f^{<k+1}(\rho_{\ker(\lambda_k),(W,\phi)})$ and $f^{<k+1}(\rho_{\coker(\lambda_k),(W,\phi)})$ are isomorphisms, then $f^{<k+1}(\rho_{\ker(\lambda_k),(W,\phi)})$, $f^{<k+1}(\rho_{\coker(l_k(K),(W,\phi)})$ and\\ $f^{<k+1}(\rho_{\coker(\lambda_k),(W,\phi)})$ are all isomorphisms. Therefore, by the five lemma, $f^{<k+1}(\rho_{K,(W,\phi)})$ is also.
\end{proof}

Since $\lambda_k\ :\ l_{k+1}(K)\rightarrow l_k(K)$ is an isomorphism in $\U/\Nil_k$, both its kernel and cokernel are in $\Nil_k$. In the following section, we will give a way to compute $\textbf{C}_{k+1}(M;K)$ for $M$ a $k$-nilpotent object in $K-\U$. Also, the authors of \cite{HLS2} proved that, for $K$ noetherian and $t$ big enough $\lambda_t\ :\ K\overset{\cong}{\rightarrow} l_t(K)$. Hence, for $K$ noetherian, Theorem \ref{gaga} provides us with an algorithm that computes $\textbf{C}(K)$, under the condition that we can compute $\textbf{C}_1(K)$ as well as $\textbf{C}_{k+1}(M;K)$ for any $k\leq t-1$ and any $k$-nilpotent object $M$ in $K-\U$.\\

For $M\in K-\U$, we can make a similar construction. We consider the following exact sequence in $K-\U$: $$0\rightarrow\ker(\lambda_k)\rightarrow l_{k+1}(M)\overset{\lambda_k}{\rightarrow}l_k(M)\rightarrow \coker(\lambda_k)\rightarrow 0$$ and get the following theorem whose proof is a straightforward adaptation of the proof of Theorem \ref{gaga}.

\begin{theorem}\label{gaga2}
Let $K$ be an unstable algebra, $\phi\in\HomK(K,H^*(W))$ and $M\in K-\U$. Then, the two following conditions are equivalent:
\begin{enumerate}
    \item $(W,\phi)\in\textbf{C}_{k+1}(M;K)$,
    \item $(W,\phi)\in\textbf{C}_k(M;K)$, $(W,\phi)\in\textbf{C}_{k+1}(\ker(\lambda_k);K)$ and  $(W,\phi)\in\textbf{C}_{k+1}(\coker(\lambda_k);K)$.
\end{enumerate}
\end{theorem}

\subsection{Central elements of functors}

In order to compute the centre of an unstable algebra using Theorem \ref{gaga}, we need to be able to compute $\textbf{C}_1(K)$ as well as $\textbf{C}_{k+1}(M;K)$ for $M$ a $k+1$-nilpotent object in $K-\U$. In this section, we define a notion of central element for objects in two different functor categories. This gives an efficient way to compute those.\\ 

We start by recalling the definition of central elements for objects in $\PFVF$, already defined and studied in \cite{B1}, as well as their relationship (also studied in \cite{B1}) with $\textbf{C}_1(K)$.

\begin{definition}
For $F\in\PFVF$ and $W\in\E$, let $\Sigma_W F\in\PFVF$ be the functor which maps $V\in\E$ to $F(W\oplus V)$ and $\alpha\in\Hom_{\E}(U,V)$ to $F(\id_W\oplus\alpha)$ from $\Sigma_WF(V)$ to $\Sigma_WF(U)$.
\end{definition}

\begin{remark}
For $F\in\PFVF$ and $W\in\E$, $\F^{\Sigma_W F}\cong \Delta_W \F^F$. In particular, for $K\in\K$, $f(T_W(K))\cong\F^{\Sigma_W g(K)}$, hence $g(T_W(K))\cong \Sigma_Wg(K)$.
\end{remark}

We will now characterise, in term of the functor $g(K)$, the fact that $f^{<1}(\rho_{K,(W,\phi)})$ is an isomorphism.

\begin{definition}
For $F\in\PFVF$ and $W$ and $V\in\E$, the injection from $0$ to $V$ induces a surjective morphism of profinite sets, natural in $W$, from $\Sigma_W F(V)$ to $F(W)$.
\begin{enumerate}
    \item For $\phi\in F(W)$, let $\Sigma_{(W,\phi)} F (V)$ be the fibre of this surjection over $\{\phi\}$.
    \item Let also $\rho_{F,W}\ : \ \Sigma_W F(V)\rightarrow \Sigma_0 F(V)\cong F(V)$ be the natural morphism induced by the injection from $0$ to $W$.
    \item Finally, let $\rho_{F,(W,\phi)}$ be the composition of $\rho_{F,W}$ with the injection from $\Sigma_{(W,\phi)}F$ to $\Sigma_{W}F$.
\end{enumerate}

\end{definition}

\begin{definition}
For $F\in\PFVF$, define $\textbf{C}(F)$ to be the set of pairs $(W,\phi)$ with $\phi\in F(W)$ such that $\rho_{F,(W,\phi)}$ is an isomorphism.
\end{definition}

The projection from $T_W(K)$ to $T_{(W,\phi)}(K)$ induces an injection from $\HomK(T_{(W,\phi)}(K),H^*(V))$ to $\HomK(T_W(K),H^*(V))\cong \Sigma_W \HomK(K,H^*(V))$. We identify this injection in terms of the functor $g(K)$.

\begin{lemma}
\begin{enumerate}
    \item For $K\in\K$, $V$ and $W$ in $\E$ and $\phi\in\HomK(K,H^*(W))$, $\HomK(T_{(W,\phi)}(K),H^*(V))$ identifies with $\Sigma_{(W,\phi)}\HomK(K,H^*(V))$ as a sub-set of $\HomK(T_W(K),H^*(V))\cong\Sigma_W\HomK(K,H^*(V))$.
    \item Under this identification, and for $F:=\HomK(K,H^*(\_))$, $$g(\rho_{K,(W,\phi)})\ :\ \HomK(T_{(W,\phi)}(K),H^*(\_))\rightarrow \HomK(K,H^*(\_))$$ identifies with $\rho_{F,(W,\phi)}$.
\end{enumerate}

\end{lemma}

\begin{proof}
The image of the injection from $\HomK(T_{(W,\phi)}(K),H^*(V))$ to $\HomK(T_W(K),H^*(V))$ is the set of  morphisms from $T_W(K)$ to $H^*(V)$, which factorise through $T_W(K)\twoheadrightarrow T_{(W,\phi)}(K)$. For $\alpha$ from $T_W(K)$ to $H^*(V)$, this is the case if and only if the composition $$T_W(K)\rightarrow H^*(V)\overset{\epsilon_{H^*(V)}}{\longrightarrow} \F,$$ is the adjoint of $\phi$. By adjunction, this is the same as asking that, for $$\Tilde{\alpha}\in\Sigma_W\HomK(K,H^*(V))\cong\HomK(K,H^*(W)\otimes H^*(V))$$ the adjoint of $\alpha$, the composition $$K\overset{\Tilde{\alpha}}{\rightarrow} H^*(W)\otimes H^*(V)\overset{\id\otimes\epsilon_{H^*(V)}}{\longrightarrow}H^*(W)\otimes\F,$$ is equal to $\phi$. Up to the isomorphism $\HomK(K,H^*(W)\otimes H^*(V))\cong\Sigma_W\HomK(K,H^*(V))$, $\alpha$ factorises through $T_W(K)\twoheadrightarrow T_{(W,\phi)}(K)$, if and only if $\Tilde{\alpha}\in\Sigma_{(W,\phi)}\HomK(K,H^*(V))$.\\

The second point is straightforward.

\end{proof}

$\rho_{g(K),(W,\phi)}=g(\rho_{K,(W,\phi)})\ :\ g(T_{(W,\phi)}(K))\rightarrow g(K)$ is therefore an isomorphism if and only if for every $\psi\ :\ K\rightarrow H^*(V)$, there is a unique morphism of unstable algebras $\phi\boxplus\psi\ :\ K\rightarrow H^*(W)\otimes H^*(V)$, such that the composition of $\phi\boxplus\psi$ with $\id\otimes\epsilon_{H^*(V)}$ (respectively $\epsilon_{H^*(W)}\otimes \id$) is equal to $\phi$ (respectively $\psi$). We get the following proposition as a consequence.

\begin{proposition}
For $K\in\K$, $\textbf{C}_1(K)=\textbf{C}(g(K))$. 
\end{proposition}

We give a similar description of $\textbf{C}_{k+1}(M;K)$ in terms of the functor $\Hom_{K-\U}(M,I_{(\_,\_)}(k))\in(\E)^{\mathcal{S}(K)^{\text{op}}}$, when $M$ is $k$-nilpotent and when $\HomK(K,H^*(\_))$ takes values in the category of finite sets. 

By adjunction and by Proposition \ref{keykey}, $f^k(T_{(W,\phi)}(M))\cong\bigoplus\limits_{\psi\in\Sigma_{(W,\phi)}g(K)(V)}\Hom_{K-\U}(M,I_{(W\oplus V,\psi)}(k))^\sharp$. Under this isomorphism, $f^k(\rho_{M,(W,\phi)})$ is nothing but the dual of the direct sum of the applications from $\Hom_{K-\U}(M,I_{(W\oplus V,\psi)}(k))$ to $\Hom_{K-\U}(M,I_{(V,(\iota_V^{W\oplus V})^*\psi)}(k))$ induced by $\iota_V^{W\oplus V}$. We get from this observation, the appropriate notion of central element for an object of $\Skop{K}$.

\begin{definition}\label{defc3}
For $K$ an unstable algebra such that $\HomK(K,H^*(\_))$ takes values in finite sets and $F\in\Skop{K}$, we say that an element $(W,\phi)\in\mathcal{S}(K)$ is $F$-central if, for any $x\in F(V,\psi)$, there exists a unique family $(x_i,\psi_i)_{1\leq i\leq n}$ such that \begin{enumerate}
    \item $\psi_i\in\Sigma_{(W,\phi)}\HomK(K,H^*(V))$,
    \item $(\iota_V^{W\oplus V})^*\psi_i=\psi$,
    \item $x_i\in F(W\oplus V,\psi_i)$,
    \item $x=\sum\limits_{1\leq i\leq n}(\iota_V^{W\oplus V})^*x_i$.
\end{enumerate}

We denote by $\textbf{C}(F;K)$ the set of central elements of $F$.
\end{definition}

\begin{proposition}\label{oopoop}
For $K$ noetherian and $M\in K-\Nil_k$ finitely-generated as a $K$-module, $\textbf{C}_{k+1}(M;K)=\textbf{C}(\Hom_{K-\U}(M,I_{(\_,\_)}(k));K)$.
\end{proposition}

\begin{proof}
This is a straightforward consequence of the discussion before Definition \ref{defc3} and Proposition \ref{keykey}.
\end{proof}

For $(W,\phi)\in\textbf{C}_1(K)$, the centrality condition for $F\in\Skop{K}$ becomes much simpler.

\begin{proposition}\label{remimp}
For $(W,\phi)\in\textbf{C}_1(K)$ and $F\in\Skop{K}$, $(W,\phi)$ is central for $F$, if for $x\in F(V,\psi)$ there is a unique $\Tilde{x}\in F(W\oplus V,\phi\boxplus\psi)$ such that $(\iota_V^{W\oplus V})^*\Tilde{x}=x$.
\end{proposition}

\begin{proof}
We recall that $\phi\boxplus\psi$ is the unique element in $\Sigma_{(W,\phi)}\HomK(K,H^*(V))$ such that $(\iota_V^{W\oplus V})^*\phi\boxplus\psi=\psi$, the proof is then straightforward.
 
\end{proof} 

For $F=\Hom_{K-\U}(M,I_{(\_,\_)}(k))$, this means that $(W,\phi)$ is $F$-central if for every morphism of $K$-modules $\alpha$ from $M$ to $I_{(V,\psi)}(k)$, there exists a unique morphism $\Tilde{\alpha}$ in $K-\U$ such that the following diagram commutes : 
$$\xymatrix{ & I_{(W\oplus V,\phi\boxplus\psi)}(k)\ar[dd]^{(\iota_V^{W\oplus V})^*}\\
M\ar[ru]^{\Tilde{\alpha}}\ar[rd]_\alpha &\\
& I_{(V,\psi)}(k).}$$

\section{Examples}

In this section we give examples of the computation of $\textbf{C}_1(M;K)$ as well as examples of the computations of $\textbf{C}(K)$ using the nilpotent filtration.

\subsection{Examples of computations of $\textbf{C}_1(M;K)$}

In this subsection, we give several examples of computations of $\textbf{C}_1(M;K)$ for some choice of $K\in\K$ $nil_1$-closed and $M\in K-\U$. In \cite{B1}, we defined the $nil_1$-closed, integral, noetherian, unstable algebra $H^*(U)^\mathcal{G}$, for $\mathcal{G}$ an appropriate groupoid with objects the sub-vector spaces of $U$. This definition is a generalisation of the algebra of invariant $H^*(U)^G$, when $G$ is a group, and we proved that any $nil_1$-closed, integral, noetherian, unstable algebra is isomorphic to some $H^*(U)^\mathcal{G}$. Here, we only give examples where $K\cong H^*(U)$ or $K\cong H^*(U)^G$.\\ 

We start by giving the example of $\textbf{C}_1(M;H^*(U))$ for $M$ a sub-object of $H^*(V,\gamma^*)$, for $\gamma$ a morphism from $V$ to $U$. We will begin by computing $\textbf{C}_1(H^*(V,\gamma^*);H^*(U))$ (which is equal to $\textbf{C}(H^*(V,\gamma^*);H^*(U))$ since $H^*(V)$ is $nil_1$-closed) and we follow with its submodules.\\

Since, $\textbf{C}(H^*(U))=\mathcal{S}(H^*(U))$, we can use the criterion for centrality from Proposition \ref{remimp}.  

\begin{proposition}\label{exa1}
For $U$, $V$ and $W\neq 0$ vector spaces and for $\gamma$ and $\psi$ linear maps from $V$ and $W$ to $U$, $(W,\psi^*)\in\textbf{C}_1(H^*(V,\gamma^*);H^*(U))$ if and only if
\begin{enumerate}
    \item $\gamma$ is injective,
    \item $\Ima(\psi)\subset\Ima(\gamma)$.
\end{enumerate}
\end{proposition}

\begin{proof}
Indeed, since, $H^*(V,\gamma^*)$ is an unstable $H^*(U)$-algebra, we can use the criterion from Proposition \ref{youpcentr}. Since $\textbf{C}(H^*(V))=\mathcal{S}(H^*(V))$, $(W,\psi^*)$ is $H^*(V,\gamma^*)$-central if and only if there exists a unique $\phi$ from $W$ to $V$ such that $\phi^*\circ\gamma^*=\psi^*$, which is the case if and only if $\gamma\circ\phi=\psi$. The condition that $\Ima(\psi)\subset \Ima(\gamma)$ is necessary and sufficient for $\phi$ to exist; injectivity of $\gamma$ is equivalent to the unicity of $\phi$.
\end{proof}

In Proposition \ref{exa1} we can interpret the two conditions in terms of the functor\\ $\Hom_{H^*(U)-\U}(H^*(V,\gamma)^*,H^*(\_,\_))$.  We consider the identity of $H^*(V,\gamma^*)$, and $(W,\psi^*)\in\mathcal{S}(U)$ such that $\Ima(\psi)\not\subset\Ima(\gamma)$. In this case, for $u\in U^\sharp$, such that $u|_{\Ima(\gamma)}=0$ and $u|_{\Ima(\psi)}\neq 0$, $\gamma^*u=0$ and $\gamma^*\boxplus \psi^*u\neq 0$, hence there can be non non-trivial morphism of unstable $H^*(U)$-module from $H^*(V,\gamma^*)$ to $H^*(V\oplus W,\gamma^*\boxplus\psi^*)$.\\

Alternatively, if the first condition does not hold, we consider $S$ a sub-vector space of $V$ such that $V=\ker(\gamma)\oplus S$. Then, for any $\phi$ from $W$ to $S$ such that $\gamma\circ\phi=\psi$ and for any $\alpha$ from $W$ to $\ker(\gamma)$, the following diagram of unstable $H^*(U)$-modules commutes: $$\xymatrix{ & H^*(W,\psi^*)\ar[dd]^{0^*} \\
H^*(V,\gamma^*)\ar[ru]^{(\phi\oplus\alpha)^*}\ar[rd]_{0^*} & \\
& H^*(0,0^*).}$$ Unless $W=0$, this negates the unicity of $\phi^*$ from $H^*(V,\gamma^*)$ to $H^*(0\oplus W,0^*\boxplus\psi^*)$ such that the following diagram commutes: 
$$\xymatrix{ & H^*(0\oplus W,0^*\boxplus\psi^*)\ar[dd]^{(\iota_0^{0\oplus W})^*} \\
H^*(V,\gamma^*)\ar[ru]^{\phi^*}\ar[rd]_{0^*} & \\
& H^*(0,0^*).}$$

For $M$ a sub-module of $H^*(V,\gamma^*)$, since $H^*(W,\psi^*)$ is injective (Proposition \ref{injmodK}), any morphism in $H^*(U)-\U$ from $M$ to $H^*(W,\psi^*)$ can be extended into a morphism from $H^*(V,\gamma^*)$ to $H^*(W,\psi^*)$. Hence, the existence condition in Proposition \ref{remimp} is satisfied for $\Hom_{H^*(U)}(M,H^*(\_,\_))$ if and only if it is satisfied for $\Hom_{H^*(U)}(H^*(V,\gamma^*),H^*(\_,\_))$. In Proposition \ref{exa2}, we prove that the uniqueness condition is satisfied, for $M\subset\gamma^*(H^*(U))$.\\ 

Conversely, in Proposition \ref{exa8}, we prove that for $M\not\subset\gamma^*(H^*(U))$ the different morphisms $(\phi\oplus\alpha)^*$ of the precedent discussion define different morphisms when restricted to $M$ and that therefore the uniqueness condition of Proposition \ref{remimp} can never be satisfied.

\begin{proposition}\label{exa2}
Let $M\in H^*(U)-\U$ a sub-module of $H^*(V,\gamma^*)$, for $U$, $V$ and $W\neq 0$ vector spaces and $\gamma$ and $\psi$ linear maps from $V$ and $W$ to $U$. Suppose that $M\subset\gamma^*(H^*(U))$. Then, $(W,\psi^*)\in\textbf{C}_1(M;H^*(U))$ if and only if $\Ima(\psi)\subset\Ima(\gamma)$.
\end{proposition}

\begin{proof}
Since $H^*(W,\psi^*)$ is injective in $H^*(U)-\U$, the existence condition from Definition \ref{defc3} is true for $\Hom_{H^*(U)-\U}(M,H^*(\_,\_))$ if and only if it is true for $\Hom_{H^*(U)-\U}(H^*(V,\gamma^*),H^*(\_,\_))$. \\ 

Let us show that the unicity is always verified if $M\subset \gamma^*(H^*(U))$. We consider a morphism $\alpha$ in $H^*(U)-\U$ from $M$ to some $H^*(X,\xi^*)$. If there exists $\Tilde{\alpha}$ from $M$ to $H^*(X\oplus W,\xi^*\boxplus\psi^*)$ such that the following diagram commutes : $$\xymatrix{ & H^*(X\oplus W,\xi^*\boxplus\psi^*)\ar[dd]^{(\iota_X^{X\oplus W})^*}\\
M\ar[ru]^{\Tilde{\alpha}}\ar[rd]_\alpha &\\
& H^*(X,\xi^*),}$$ since the $H^*(Y,\upsilon^*)$ are injective, there exist also $\beta$ and $\Tilde{\beta}$ whose restrictions to $M$ are $\alpha$ and $\Tilde{\alpha}$ and such that the following diagram commutes : $$\xymatrix{ & H^*(X\oplus W,\xi^*\boxplus\psi^*)\ar[dd]^{(\iota_X^{X\oplus W})^*}\\
H^*(V,\gamma^*)\ar[ru]^{\Tilde{\beta}}\ar[rd]_\beta &\\
& H^*(X,\xi^*).}$$ The choice of $\beta$ and $\Tilde{\beta}$ need not be unique, but since $M\subset \gamma^*(H^*(U))$ by hypothesis, and since the image of $\gamma^*(H^*(U))$ by any morphism in $H^*(U)-\U$ from $H^*(V,\gamma^*)$ is determined from the image of $1_{H^*(V)}$, the restriction to $M$ of any such $\Tilde{\beta}$ will coincide with $\Tilde{\alpha}$. Hence, $\Tilde{\alpha}$ is unique.
\end{proof}

\begin{proposition}\label{exa8}
Let $M\in H^*(U)-\U$ be a sub-module of $H^*(V,\gamma)^*$, $U$, $V$ and $W\neq 0$ vector spaces and $\gamma$ and $\psi$ linear maps from $V$ and $W$ to $U$. If $M\not\subset\gamma^*(H^*(U))$, $(W,\psi^*)\not\in\textbf{C}_1(M;H^*(U))$.
\end{proposition}

\begin{proof}

If $\Ima(\psi)\not\subset\Ima(\gamma)$, $(W,\psi^*)$ is not central for the same reason as in Proposition \ref{exa2}. We consider the inclusion $M\hookrightarrow H^*(V,\gamma^*)$. Then, for $S$ such that $S\oplus \ker(\gamma)=V$, we consider $\phi\ :\ W\rightarrow V$ the unique map such that $\Ima(\phi)\subset S$ and $\gamma\circ\phi=\psi$. Then, for every $\alpha\ :\ W\rightarrow V$ such that $\Ima(\alpha)\subset\ker(\gamma)$, the morphism $\phi_\alpha\ :\  W\oplus V \rightarrow V$ which maps $x\oplus y\in W\oplus V$ to $y+\phi(x)+\alpha(x)\in V$ satisfies: \begin{enumerate}
    \item $\phi_\alpha^*\ :\ H^*(V,\gamma^*)\rightarrow H^*(W\oplus V,\psi^*\boxplus\gamma^*)$ is a morphism of $H^*(U)$-modules, 
    \item the following diagram commutes : $$\xymatrix{ & H^*(W\oplus V,\psi^*\boxplus\gamma^*)\ar[dd]^{(\iota_V^{W\oplus V})^*}\\
M\ar[ru]^{\phi_\alpha^*}\ar@{^(->}[rd] &\\
& H^*(V,\gamma^*).}$$
\end{enumerate}

Let us show that no element in $H^*(V,\gamma^*)\backslash \gamma^*(H^*(U))$ has the same image under each $\phi_\alpha^*$. We consider the decomposition $$H^*(V)\cong H^*(S)\otimes H^*(\ker(\gamma))\cong H^*(S)\oplus\Bar{H}^*(\ker(\gamma))\otimes H^*(S),$$ where $\Bar{H}^*(\_)$ denote the reduced cohomology over $\F$. Any element in $x\in H^*(V)$ can be written in a unique way as $x=s\oplus(\bigoplus_{b\in\mathcal{B}}r_b\otimes b)$ with $s\in H^*(S)$, $\mathcal{B}$ a basis of $H^*(S)$ and $r_b\in\Bar{H}^*(\ker(\gamma))$. Therefore, for $\alpha$ and $\beta$ morphisms from $W$ to $\ker(\gamma)$, $(\phi_\alpha^*-\phi_\beta^*)(x)=\bigoplus_{b\in\mathcal{B}}(\phi_\alpha^*-\phi_\beta^*)(r_b)\otimes \phi_\alpha^*b$. This is because, the restriction of $\phi_\alpha^*$ to $H^*(S)$ does not depend on $\alpha$. Hence, $(\phi_\alpha^*-\phi_\beta^*)(x)=0$ if and only if $(\phi_\alpha^*-\phi_\beta^*)(r_b)=0$ for each $b\in\mathcal{B}$. So we only have to prove that there is no element in $x\in\Bar{H}^*(\ker(\gamma))$ such that $\phi_0^*(x)=\phi_\alpha^*(x)$ for every $\alpha$.\\

The restriction of $\phi_0^*$ to $H^*(\ker(\gamma))$ is the inclusion of $H^*(\ker(\gamma))$ in $H^*(W\oplus V)$, we consider $(v_1,...,v_n)$ a basis of the dual of $\ker(\gamma)$. For a given $i$, there is a unique way to write $x$ as $x=a_0+a_1v_i+...+a_nv_i^n$ with $n\in\N$, $(a_i)_{i\in\left[|1,n|\right]}\in \F\left[v_1,...,\hat{v}_i,...,v_n\right]^n$ and $a_n\neq 0$. And for $x\neq 1_{H^*(\ker(\gamma))}\in H^*(\ker(\gamma))$, we can chose $i$ such that $n\neq0$.  
Then, for $w$ a non trivial linear form on $W$, we can chose $\alpha$ such that $\phi_\alpha^*v_i=v_i+w$. Then, $(\phi_0^*-\phi_\alpha^*)(x)=a_1w+a_2(v_i^2-(v_i+w)^2)+...a_n(v_i^n-(v_i+w)^n)$ whose degree as a polynomial in $w$ is $n$, so $(\phi_0^*-\phi_\alpha^*)(x)\neq 0$. This concludes the proof.

\end{proof}

We now give some similar example from  $H^*(U)^G-\U$ for $G$ some sub-group of $\text{Gl}(U)$. In this case, we have that $\HomK(H^*(U)^G,H^*(V))\cong \Hom(V,U)/G$, where $G$ acts on $\Hom(V,U)$ by composition. It is worth recalling from Example \ref{gluglu}, that for $\psi$ from $W$ to $U$, and for $\left[\psi\right]$ the class of $\psi$ under the action of $G$, $(W,\left[\psi\right])\in\textbf{C}(H^*(U)^G)$ if and only if for all $y\in\Ima(\psi)$ and for all $g\in G$, $gy=y$. 

\begin{proposition}\label{exa4}
For $U$, $V$ and $W\neq 0$ some vector spaces and for $\gamma$ and $\psi$ some morphisms from $V$ and $W$ to $U$, $(W,\left[\psi\right])\in\textbf{C}_1(H^*(V,\left[\gamma\right]);H^*(U))$ if and only if
\begin{enumerate}
    \item $\gamma$ is injective,
    \item there exists $g\in G$ such that $\Ima(g\circ\psi)\subset \Ima(\gamma)$,
    \item for all $g'\in G$ such that $\Ima(g'\circ\psi)\subset\Ima(\gamma)$, $g'\circ\psi=g\circ\psi$.
\end{enumerate}
\end{proposition}

\begin{proof}
$H^*(V,\left[\gamma\right])$ is an unstable $H^*(U)^G$-algebra, so we can use the criterion from Proposition \ref{youpcentr}. Since $\textbf{C}(H^*(V))=\mathcal{S}(H^*(V))$, $(W,\left[\psi\right])$ is $H^*(V,\left[\gamma\right])$-central if and only if there exists a unique $\phi$ from $W$ to $V$ such that $\phi^*\circ\left[\gamma\right]=\left[\psi\right]$, which is the case if and only if $\gamma\circ\phi=g\circ\psi$ for some $g\in G$. The condition that $\Ima(\psi)\subset \Ima(\gamma)$ is necessary and sufficient for $\phi$ to exist, and the two other conditions ensure that $\phi$ is unique.
\end{proof}

\begin{remark}
If $(W,\psi^*)\in\textbf{C}_1(H^*(U)^G)$, then $gy=y$ for all $g\in G$ and $y\in\Ima(\psi)$. Hence, the third condition of Proposition \ref{exa4} is automatically verified.
\end{remark}

\subsection{Examples of computations of $\textbf{C}(K)$ using the nilpotent filtration}

In this subsection, we use Theorem \ref{gaga} to compute $\textbf{C}(K)$ for some examples of non $nil_1$-closed algebras. We will mainly consider algebras that are $nil_2$-closed, since the step from $\textbf{C}_1(K)=\textbf{C}(l_1(K))$ to $\textbf{C}_2(K)=\textbf{C}(l_2(K))$ is not different from the step from $\textbf{C}_2(K)$ to $\textbf{C}_3(K)$ and so on. Also, we will relate those examples to the existence of $H^*(W)$-comodule structure on $K$ from Corollary \ref{DWCT}.\\ 

In the following, $V_2$ and $V_3$ will denote the vector spaces $\F x\oplus\F y$ and $\F x\oplus\F y\oplus\F z$, $\pi$ will denote the canonical projection from $V_3$ onto $V_2$ and $B_2$ will denote the subgroup of $\text{Gl}(V_2)$ generated by the morphism which sends $x$ to itself and $y$ to $x+y$.

\begin{example}
We consider the unstable algebra $K:=H^*(V_2)^{B_2}\oplus\Sigma H^*(V_3,\left[\pi\right])$ where the product of two elements from $\Sigma H^*(V_3,\left[\pi\right])$ is $0$ and the product of $x\in H^*(V_2)^{B_2}$ by $y\in\Sigma H^*(V_3,\left[\pi\right])$ is given by the $H^*(V_2)^{B_2}$-module structure $\left[\pi\right]$ on $\Sigma H^*(V_3)$. Since $\Sigma H^*(V_3,\left[\pi\right])$ is nilpotent, the projection from $K$ to $H^*(V_2)^{B_2}$ (which is $nil_1$-closed, eg \cite{HLS2}) is an isomorphism in $\U/\Nil_1$. Hence, this is the $nil_1$-localisation of $K$. We therefore have that $\mathcal{S}(K)\cong\mathcal{S}(H^*(V_2)^{B_2})$ and for $\psi$ a linear map from a vector space $W$ to $V_2$, $(W,\left[\psi\right])\in\mathcal{S}(K)$ is in $\textbf{C}_1(K)$ if and only if $\Ima(\psi)\subset\F x$.\\

Then, by Theorem \ref{gaga}, we have that for $(W,\left[\psi\right])$ such that $\Ima(\psi)\subset\F x$, $(W,\left[\psi\right])$ is in $\textbf{C}_2(K)$ if and only if it is in $\textbf{C}_2(\Sigma H^*(V_3,\left[\pi\right]);K)$. By Proposition \ref{oopoop}, $\textbf{C}_2(\Sigma H^*(V_3,\left[\pi\right]);K)=\textbf{C}(\Hom_{K-\U}(\Sigma H^*(V_3,\left[\pi\right]),I_{(\_,\_)}(1));K)$. But using the adjunction between the functors $\Sigma$ and $\Omega$, $\Hom_{K-\U}(\Sigma H^*(V_3,\left[\pi\right]),I_{(\_,\_)}(1))=\Hom_{H^*(V_2)^{B_2}-\U}(H^*(V_3,\left[\pi\right]),H^*(\_,\_))$. Hence, by Proposition \ref{exa4}, since $\pi$ is not injective, $(W,\left[\psi\right])\in\textbf{C}_2(\Sigma H^*(V_3,\left[\pi\right]);K)=\textbf{C}_1(H^*(V_3,\left[\pi\right]);H^*(V_2)^{B_2})$ if and only if $W=0$. We find that $\textbf{C}_2(K)=\{(0,\epsilon_K)\}$. Since, $K$ is $nil_2$-closed, $\textbf{C}(K)=\textbf{C}_2(K)$.\\

It is worth noticing that, nevertheless, there is a $H^*(\F)$-comodule structure in $\K$ on $K$, such that the composition $K\rightarrow K\otimes H^*(\F)\rightarrow H^*(\F)$ is equal to the morphism which sends $\Sigma H^*(V_3)$ to $0$ and whose restriction to $H^*(V_2)^{B_2}$ is equal to $\left[\iota_x\right]$, with $\iota_x$ the inclusion from $\F x$ to $V_2$. For $(u,v,w)$ the dual basis of $(x,y,z)$, we have $K=\F\left[v,u(u+v)\right]\oplus\Sigma \F\left[u,v,w\right]$. The comodule structure is the one that send $v$ to $v\otimes 1$, $u(v+u)$ to $u(v+u)\otimes 1+v\otimes u +1\otimes u^2$, $\sigma v$ to $\sigma v\otimes 1$, $\sigma w$ to $\sigma w\otimes 1$ and, finally, $\sigma u$ to $\sigma u\otimes 1+\sigma 1\otimes u$, where $\sigma x\in\Sigma \F\left[u,v,w\right]$ denotes the suspension of $x\in \F\left[u,v,w\right]$. There is no contradiction with Corollary \ref{DWCT} since, in this case, $\mathcal{Q}(K)$ is not locally finite, indeed the projection from $I(K)$ to $\mathcal{Q}(K)$ induces an injection from  $\A\sigma w$ to $\mathcal{Q}(K)$ and $\A\sigma w$ is not finite.
\end{example}

The two following examples are similar to the first one, when we replace $\Sigma H^*(V_3,\left[\pi\right])$ by other modules of the form $\Sigma H^*(V,\left[\gamma\right])$, in order to illustrate the relevance of the other criteria of Proposition \ref{exa4}.

\begin{example}
This time, we consider the unstable algebra $K:=H^*(V_2)^{B_2}\oplus\Sigma H^*(V_2,\left[\id_{V_2}\right])$. We again have that the $nil_1$-localisation of $K$ is the projection onto $H^*(V_2)^{B_2}$. Thus $\mathcal{S}(K)\cong\mathcal{S}(H^*(V_2)^{B_2})$ and, for $\psi$ a linear map from a vector space $W$ to $V_2$, $(W,\left[\psi\right])\in\mathcal{S}(K)$ is in $\textbf{C}_1(K)$ if and only if $\Ima(\psi)\subset\F x$.\\

We then compute $\textbf{C}_2(K)$ using Theorem \ref{gaga}; as in the previous example, we have that $(W,\psi^*)\in\textbf{C}_1(K)$ is in $\textbf{C}_2(K)$ if and only if it is in $\textbf{C}_1(H^*(V_2,\left[\id_{V_2}\right]);K)$. But this time, the three conditions of Proposition \ref{exa4} are always satisfied by pairs $(W,\left[\psi\right])$ in $\textbf{C}_1(K)$, therefore, since $K$ is $nil_2$-closed, $\textbf{C}(K)=\textbf{C}_1(K)$.\\

In this case, $\mathcal{Q}(K)$ is finite (in particular it is locally finite), and the only non trivial $H^*(\F)$-comodule structure on $H^*(V_2)^{B_2}$, which sends $v$ to $v\otimes 1$ and $u(v+u)$ to $u(v+u)\otimes 1+v\otimes u+1\otimes u^2$, can be extended to $K\cong \F\left[v,u(u+v)\right]\oplus \Sigma\F\left[u,v\right]$ by sending $\sigma v$ to $\sigma v\otimes 1$ and $\sigma u$ to $\sigma u\otimes 1+\sigma 1\otimes u$.
\end{example}

\begin{example}
Finally, we consider the unstable algebra $K:=H^*(V_2)^{B_2}\oplus\Sigma H^*(\F y,\left[\iota_y\right])$. We again have that the $nil_1$-localisation of $K$ is the projection onto $H^*(V_2)^{B_2}$ and $\mathcal{S}(K)\cong\mathcal{S}(H^*(V_2)^{B_2})$. For $\psi$ a linear map from a vector space $W$ to $V_2$, $(W,\left[\psi\right])\in\mathcal{S}(K)$ is in $\textbf{C}_1(K)$ if and only if $\Ima(\psi)\subset\F x$.\\
As above, we compute $\textbf{C}_2(K)$ using Theorem \ref{gaga}, and like in the two previous examples, we have that $(W,\psi^*)\in\textbf{C}_1(K)$ is in $\textbf{C}_2(K)$ if and only if it is in $\textbf{C}_1(H^*(\F y,\left[\iota_y\right]);K)$. This time, the second condition of Proposition \ref{exa4} is not satisfied for $(W,\left[\psi\right])\in\textbf{C}_1(K)$ with $\psi$ non trivial, indeed the second condition of Proposition \ref{exa4} would require $\Ima(\psi)$ to be a sub-vector space of $\F y$, when $(W,\left[\psi\right])\in\textbf{C}_1(K)$ implies that $\Ima(\psi)\subset \F x$. Hence, $\textbf{C}(K)=\{(W,\epsilon_{K,W})\ ;\ W\in\E\}$.\\

In this case, $\mathcal{Q}(K)$ is locally finite, and the non trivial $H^*(\F)$-comodule structure on $H^*(V_2)^{B_2}$, cannot be extended to $K\cong \F\left[v,u(u+v)\right]\oplus \Sigma\F\left[v\right]$. Indeed, $\sigma v$ would need to be sent to $\sigma v\otimes 1$, hence $u(v+u)\sigma v=0$ would require to be sent to $\sigma v^2\otimes u+\sigma v\otimes u^2\neq 0$.
\end{example}

In the above three example, we computed $\textbf{C}_2(K)$ from $\textbf{C}_1(K)$ and $\textbf{C}_2(nil_1(K)/nil_2(K);K)$. Let us give an example when the obstruction to centrality comes from $\coker(\lambda_1)$.

\begin{example}
We consider $K=\F\left[v,w\right]\oplus u^2\F\left[u,v,w\right]$. The injection from $K$ to $H^*(V_3)\cong\F\left[u,v,w\right]$ has cokernel $\Sigma\F\left[v,w\right]$ which is nilpotent. Therefore, the former injection is the $nil_1$-localisation of $K$. We therefore have that $\mathcal{S}(K)\cong\mathcal{S}(H^*(V_3))$ and that, under this identification, $\textbf{C}_1(K)=\mathcal{S}(H^*(V_3))$. So, for any linear map $\psi$ from a vector space $W$ to $V_3$, $(W,\psi^*\circ\lambda_1)\in\textbf{C}_1(K)$.\\

As before, $(W,\psi^*\circ\lambda_1)\in\textbf{C}_2(K)$ if and only if it is in  $\textbf{C}_2(\Sigma\F\left[v,w\right];K)$, which is the case if and only if it is in $\textbf{C}_1(\F\left[v,w\right];H^*(V_3))$. By proposition \ref{exa1}, this is the case if and only if $\Ima(\psi)\subset\F y\oplus \F z$. Since $K$ is $nil_2$-closed, by Theorem \ref{gaga} we get that $\textbf{C}(K)=\{(W,\psi^*\circ\lambda_1)\ ;\ \Ima(\psi)\subset\F y\oplus \F z\}$.\\

We can again interpret this result in terms of Corollary \ref{DWCT}. We consider the $\F\left[u\right]$-comodule structure on $H^*(V_3)$ in $\K$ that send $v$ to $v\otimes 1$, $w$ to $w\otimes 1$ and $u$ to $u\otimes 1+1\otimes u$, and which is associated to $(\F x,\iota_x)\in\textbf{C}(H^*(V_3))$. This comodule structure does not induce a comodule structure on $K$, since the image of $u^3$, which is in $K$, is $u^3\otimes 1+ u^2\otimes u+u\otimes u^2+1\otimes u^3$ which is not in $K\otimes \F\left[u\right]$.
\end{example}

\bibliographystyle{alpha}
\bibliography{Biblio.bib}

\begin{center}
    Address : Université de Lille, laboratoire Paul Painlevé (UMR8524), Cité scientifique, bât. M3, 59655
VILLENEUVE D’ASCQ CEDEX, FRANCE\\
    e-mail : aacde13@live.fr
\end{center}

\end{document}